\documentclass{article}

\usepackage{PRIMEarxiv}
\usepackage[utf8]{inputenc} 
\usepackage[T1]{fontenc}    
\usepackage{hyperref}       
\usepackage{url}            
\usepackage{booktabs}       
\usepackage{amsfonts}       
\usepackage{nicefrac}       
\usepackage{microtype}      
\usepackage{lipsum}
\usepackage{fancyhdr}       
\usepackage{graphicx}       
\graphicspath{{media/}}     
\usepackage{amsmath,amsfonts}
\usepackage{cleveref}
\usepackage{algorithmic}
\usepackage{algorithm}
\usepackage{array}
\usepackage[caption=false,font=normalsize,labelfont=sf,textfont=sf]{subfig}
\usepackage{textcomp}
\usepackage{stfloats}
\usepackage{url}
\usepackage{verbatim}
\usepackage{graphicx}
\usepackage{ntheorem}
\theoremseparator{.}
\usepackage{cite}
\usepackage{multirow}
\usepackage{float}

\newtheorem{assumption}{Assumption}

\newtheorem{theorem}{Theorem}
\newtheorem{lemma}{Lemma}
\newtheorem{definition}{Definition}
\newtheorem{remark}{Remark}
\newtheorem{proof}{Proof}
\newtheorem{corollary}{Corollary}
\title{Alternating Stochastic Variance-Reduced Algorithms with Optimal Complexity for Bilevel Optimization}

\author{
  Haimei Huo \\
 School of Mathematical Sciences, Dalian University of Technology, Dalian, China  \\ 
\texttt{\{Haimei Huo\}email@ab1234@mail.dlut.edu.cn} \\
   \And
  Zhixun Su \\
  School of Mathematical Sciences, Dalian University of Technology, Dalian, China  \\ 
  \texttt{email@zxsu@dlut.edu.cn} \\  
}

\begin{document}
\maketitle

\begin{abstract}
This paper studies the unconstrained nonconvex-strongly-convex bilevel optimization problem. A common approach to solving this problem is to alternately update the upper-level and lower-level variables using (biased) stochastic gradients or their variants, with the lower-level variable updated either one step or multiple steps. In this context, we propose two alternating stochastic variance-reduced algorithms, namely ALS-SPIDER and ALS-STORM, which introduce an auxiliary variable to estimate the hypergradient for updating the upper-level variable. ALS-SPIDER employs the SPIDER estimator for updating variables, while ALS-STORM is a modification of ALS-SPIDER designed to avoid using large batch sizes in every iteration. Theoretically, both algorithms can find an $\epsilon$-stationary point of the bilevel problem with a sample complexity of $O(\epsilon^{-1.5})$ for arbitrary constant number of lower-level variable updates. To the best of our knowledge, they are the first algorithms to achieve the optimal complexity of $O(\epsilon^{-1.5})$ when performing multiple updates on the lower-level variable. Numerical experiments are conducted to illustrate the efficiency of our algorithms. 
\end{abstract}

\section{Introduction}\label{sec:1}
Bilevel optimization has become a mainstream framework in many machine learning applications such as hyperparameter optimization\cite{14, 1, 5}, meta-learning\cite{2, 4, 3}, and reinforcement learning\cite{7, 6}. In this paper, we consider the nonconvex-strongly-convex bilevel optimization problem given by 
\begin{align}\label{eq1}
&\min\limits_{x \in \mathbb{R}^p} \Phi(x) := f\left( x, y^*(x)\right) := \mathbb{E}_{\xi }\left[F(x, y^*(x); \xi)\right] \\  \nonumber
&\text{s.t.} ~ y^*(x) :=  { \underset { y \in \mathbb{R}^q} { \operatorname {arg\,min} } \,  \{g(x, y) := \mathbb{E}_{\zeta}[G( x, y; \zeta)]}\} 
\end{align} 
where the lower level (LL) function $g(x, y)$ is strongly convex in $y$, and the upper level value function $\Phi(x):= f\left( x, y^*(x)\right)$ is nonconvex in $x$. Here, $f(x, y)$ and $g(x, y)$ take the expected values with respect to (w.r.t.) random variables $\xi$ and $\zeta$, respectively.

The objective of problem (\ref{eq1}) is to minimize $\Phi(x)$ on $x \in \mathbb{R}^p$. Under the following assumption (i.e., Assumption \ref{assum:1}), 
$\Phi(x)$ can be ensured to be differentiable through the implicit function theorem, and for any $x \in \mathbb{R}^p$, the hypergradient $\nabla \Phi(x)$ takes the form of 
\begin{equation} \label{eq2}
\nabla \Phi(x) = \nabla_{x}f(x, y^*(x)) - \nabla_{x} \nabla_{y} g(x, y^*(x))[\nabla_{y}^2 g(x, y^*(x))]^{-1} \nabla_{y}f(x, y^*(x))
\end{equation}
where $y^*(x)$ is the optimal solution of the LL problem $\min_{ y\in \mathbb{R}^q}g(x, y)$. The definitions of $\nabla_{x} f(x, y^*(x))$, $\nabla_{y} f(x, y^*(x))$, $\nabla_{x} \nabla_{y} g(x, y^*(x))$, and $\nabla_{y}^2 g(x, y^*(x))$ can be found in Section \ref{sec:2} of this paper, and the detailed proof is available in \cite[Lemma 2.1]{16}.

\begin{assumption}\label{assum:1}
The functions  $f(x, y)$, $F(x, y; \xi)$ (resp. $g(x, y)$, $G(x, y; \zeta)$), with $\xi$ (resp. $\zeta$) being arbitrary, are continuously differentiable (resp. twice continuously differentiable) on $\mathbb{R}^p \times \mathbb{R}^q$. Moreover, for any $x$ and $\zeta$, functions $g(x, y)$ and $G(x, y; \zeta)$ are $\mu$-strongly convex in $y$; function $\Phi(x)$ is nonconvex in $x$, and $\inf_{x\in \mathbb{R}^p}\Phi(x)> -\infty$.
\end{assumption}

The formulation in (\ref{eq2}) enables us to solve problem (\ref{eq1}) by applying stochastic gradient descent (SGD) updates to the variable $x$. However, estimating the hypergradient $\nabla \Phi(x)$ via (\ref{eq2}) requires the optimal solution $y^*(x)$ of the LL problem, which is typically unknown. This makes the practical implementation of the SGD method challenging. 

To address this challenge, various stochastic algorithms have been proposed that alternately update the UL variable $x$ and the LL variable $y$ in order to solve problem (\ref{eq1}). In these algorithms, the LL variable is updated using the stochastic gradient (or its variants) of the LL function to approximate $y^*(x)$, while a biased stochastic estimate of $\nabla \Phi(x)$ (or its variants) is used to update the UL variable.

Some algorithms in this direction employ one update for the LL variable. In \cite{7}, the algorithm updates the variable $y$ using one SGD iteration, yielding an approximation $\bar{y}$ of $y^*(x)$. This leads to the following surrogate for $\nabla \Phi(x)$, denoted as $\nabla \bar{\Phi}(x)$,
\begin{equation} \label{eq3}
\nabla \bar{\Phi}(x) = \nabla_{x}f(x, \bar{y}) - \nabla_{x} \nabla_{y} g(x, \bar{y})[\nabla_{y}^2 g(x, \bar{y})]^{-1} \nabla_{y}f(x, \bar{y}).
\end{equation}
The variable $x$ is then updated based on a biased stochastic estimate of $\nabla \bar{\Phi}(x)$ (which is also a biased estimate of $\nabla \Phi(x)$), obtained by utilizing the 
truncated stochastic Neumann series to approximate the Hessian inverse $[\nabla_{y}^2 g(x, \bar{y})]^{-1}$ in (\ref{eq3}). This algorithm is further refined in \cite{20, 21, 41, 42} for improving sample complexity, with momentum or variance reduction techniques incorporated into the variable updates. However, the best known complexity among them includes an additional $\log{\epsilon^{-1}}$ factor, compared to the optimal complexity of $O(\epsilon^{-1.5})$. Several alternative algorithms \cite{40, 54, 60, 61} directly compute the inverse of the variance-reduced Hessian $\nabla_{y}^2 g(x, \bar{y})$\footnote{Here and thereafter in this section, we use $\bar{y}$ to denote the approximation of $y^*(x)$ obtained from any stochastic updates applied to the variable $y$. Correspondingly, $\nabla \bar{\Phi}(x)$ in (\ref{eq3}) represents the surrogate of $\nabla \Phi(x)$, where $y^*(x)$ is replaced by $\bar{y}$.} to obtain a variance-reduced estimate of $\nabla \bar{\Phi}(x)$, which is then used to update $x$. While some of these algorithms \cite{54, 61} achieve the optimal complexity of $O(\epsilon^{-1.5})$, computing the inverse of the Hessian is nontrivial. Recently, several works \cite{34, 35, 22, 43, 61, 59} utilize one SGD iteration or its variants to approximately compute the Hessian-inverse-vector product $[\nabla_{y}^2 g(x, \bar{y})]^{-1} \nabla_{y} f(x, \bar{y})$ in (\ref{eq3}) for obtaining an estimate of $\nabla \bar{\Phi}(x)$. Among them, BSVRB$^{\text{v}2}$ \cite{61} is shown to achieve the optimal complexity of $O(\epsilon^{-1.5})$ by utilizing the STORM variance reduction technique \cite{50} for variable updates. However, this algorithm may sometimes fail to achieve good empirical performance due to the simplistic one-step update of the LL variable. Moreover, an  extension of this algorithm to perfrom multiple iterations for the LL variable is not directly.

On the other hand, some algorithms perform multiple updates on the LL variable. In \cite{16}, the authors propose a bilevel stochastic approximation (BSA) algorithm with a sample complexity of $O(\epsilon^{-3})$. BSA uses an increasing number of SGD steps to update the variable $y$, and estimates the hypergradient $\nabla \Phi(x)$ for updating the variable $x$ using the truncated stochastic Neumann series. Building upon BSA and utilizing the warm-start initialization strategy for LL updates, stocBiO \cite{13} improves the sample complexity to $O(\epsilon^{-2}\log{\epsilon^{-1}})$, with only a constant number of SGD steps required for updating $y$. A tighter analysis of stocBiO is presented in \cite{24}, allowing the algorithm to achieve the sample complexity of $O(\epsilon^{-2}\log{\epsilon^{-1}})$ without relying on large batch sizes during iterations. AmIGO \cite{32} replaces the truncated stochastic Neumann series used in stocBiO with multiple SGD iterations, resulting in a sample complexity of $O(\epsilon^{-2})$. Recently, the sample complexity is further improved to $O(\epsilon^{-1.5}\log{\epsilon^{-1}})$ in \cite{21} through the use of the truncated stochastic Neumann series for estimating $\nabla \Phi(x)$, and the incorporation of the SPIDER variance reduction technique \cite{51} for both UL and LL variable updates, under additional Lipschitz assumptions on the stochastic functions. However, this complexity result does not match that of SPIDER for nonconvex stochastic single-level optimization, where a sample complexity of $O(\epsilon^{-1.5})$ is achievable. Furthermore, the proposed algorithm VRBO \cite{21} requires a three-level loop structure, which is computationally expensive.

To this end, we aim to answer the following question:
\begin{itemize}
\item Is it possible to leverage the SPIDER variance reduction technique to design an algorithm with a two-level loop structure that achieves a complexity of $O(\epsilon^{-1.5})$ when performing multiple updates on the LL variable ?
\end{itemize}

\subsection{Main contributions}

In this paper, we provide an affirmative answer to the above question. Our main contributions are as follows. First, we propose an alternating stochastic variance-reduced algorithm using the SPIDER variance reduction technique (ALS-SPIDER), within the existing framework of alternating updates for the UL and LL variables. The key idea, inspired by recent works \cite{32, 34, 61}, is to apply multiple stochastic updates to an auxiliary variable in order to estimate the hypergradient $\nabla \Phi(x)$. This contrasts with the SPIDER-based variance-reduced algorithm VRBO \cite{21}, which utilizes the Neumann series to estimate $\nabla \Phi(x)$. In addition, unlike VRBO \cite{21}, where the variance-reduced estimate for updating the UL variable is recursively constructed along the iterative trajectory of the LL variable, ALS-SPIDER constructs this estimate from the output of the multiple LL variable updates. As a result, ALS-SPIDER exhibits a two-level loop structure.


Second, under the same assumptions and batch size setting as in \cite{21}, we demonstrate that ALS-SPIDER finds an $\epsilon$-stationary point of problem (\ref{eq1}) with a sample complexity of $O(\epsilon^{-1.5})$, regardless of the constant number of iterations used to update the LL variable \footnote{In this paper, the number of iterations for updating the LL or auxiliary variables remains fixed throughout the optimization process.}. This resolves the above question.

Third, we modify ALS-SPIDER by replacing the SPIDER estimator with the STORM estimator, thereby eliminating the need for large batch sizes in every iteration. The resultant algorithm, named ALS-STORM, is shown to achieve the same theoretical results as ALS-SPIDER, but only requires a large batch size in the first iteration.

Fourth, the convergence of ALS-SPIDER and ALS-STORM can be guaranteed for arbitrary constant number of iterations used to update the auxiliary variable, as opposed to the $O(\log{\epsilon^{-1}})$ Hessian-vector product computations required in \cite{21} for estimating $\nabla \Phi(x)$. This is the primary reason why our algorithm achieves a better complexity than the one presented in \cite{21}.

To the best of our knowledge, ALS-SPIDER and ALS-STORM are the first algorithms that achieve the optimal complexity of $O(\epsilon^{-1.5})$ when updating the LL variable multiple times. Furthermore, this is the first work to utilize the stochastic updates to estimate the hypergradient among algorithms that perform multiple variance-reduced updates for the LL variable. It is also worth noting that our algorithms are the first to achieve the complexity of $O(\epsilon^{-1.5})$ for arbitrary constant number of iterations used to solve the LL problem or estimate the hypergradient.

\section{Preliminaries}\label{sec:2}

Throughout this paper, we use $\|\cdot\|$ to represent the $L_2$ norm for vectors and spectral norm for matrices. For a function $f(x, y)$, $\nabla f(x, y)$ denotes its gradient w.r.t. $(x, y)$, while $\nabla_{x}f(x, y)$ and $\nabla_{y}f(x, y)$ refer to its partial derivatives w.r.t. $x$ and $y$, respectively. The notation $\nabla_{x}\nabla_{y}f(x, y)$(resp. $\nabla_{y}^2 f(x, y)$) refers to the Jacobian(resp. Hessian matrix) of $\nabla_{x} f(x, y)$(resp. $f(x, y)$) w.r.t. $y$. Given a closed convex set $\Omega$, we define $\text{proj}_{\Omega}(x_0) = { \operatorname {arg\,min} }_{x\in \Omega} \, \|x - x_0\|^2/2$.

Next, we make some additional assumptions for problem (\ref{eq1}).

\begin{assumption}\label{assum:2}
The functions $f(x, y)$ and $g(x, y)$ satisfy
\begin{enumerate}
\item[(a)] $f(x, y)$ and its gradient $\nabla f(x, y)$ are Lipschitz continuous with constants $C_f$ and $L_f$, respectively, w.r.t. $(x,y)$;
\item[(b)] For $g(x, y)$, the gradient $\nabla g(x, y)$, the mixed second-order partial derivative $\nabla_{x} \nabla_{y} g(x, y)$, and the second-order partial derivative $\nabla_{y}^{2}g(x, y)$ are $L_g$-, $L_{gxy}$-, and $L_{gyy}$-Lipschitz continuous w.r.t. $(x, y)$, respectively.
\end{enumerate}
\end{assumption}

\begin{assumption} \label{assum:3}
The stochastic component functions $ F(x, y; \xi)$ and $G(x, y; \zeta)$ satisfy
\begin{enumerate}
\item[(a)] The gradients $\nabla F(x, y; \xi)$, $\nabla G(x, y; \zeta)$, the mixed second-order partial derivative $\nabla_{x} \nabla_{y} G(x, y; \zeta)$, and the second-order partial derivative $\nabla_{y}^2 G(x, y; \zeta)$ are unbiased estimates of $\nabla f(x, y)$, $\nabla g(x, y)$, $\nabla_{x} \nabla_{y} g(x, y)$, and $\nabla_{y}^2 g(x, y)$, respectively;
\item[(b)] There exist constants $\sigma_f$, $\sigma_g$, $\sigma_{gxy}$, $\sigma_{gyy}$ such that 
\begin{align*}
&\mathbb{E} \| \nabla F(x, y; \xi) - \nabla f(x, y) \|^2 \le \sigma_f^2, ~~~~~\mathbb{E} \|\nabla G(x, y; \zeta) - \nabla g(x, y) \|^2 \le \sigma_g^2, \\ 
&\mathbb{E}\| \nabla_{x} \nabla_{y} G(x, y; \zeta) - \nabla_{x} \nabla_{y} g(x, y)\|^2  \le \sigma_{gxy}^2, \\
&\mathbb{E}\|\nabla_{y}^2 G(x, y; \zeta) - \nabla_{y}^2 g(x, y)\|^2  \le \sigma_{gyy}^2.
\end{align*}
\end{enumerate}
\end{assumption}

The Lipschitz continuity conditions on the objective functions in Assumption \ref{assum:2}, along with the unbiasedness and bounded variance conditions for the stochastic component functions in Assumption \ref{assum:3}, are commonly used in the analysis of nonconvex-strongly-convex bilevel optimization; see e.g., \cite{7, 13, 16, 21}. Under Assumptions \ref{assum:1}, \ref{assum:2}, we establish the Lipschitz continuity properties for $\nabla \Phi(x)$, $y^*(x)$, and $v^*(x)$ in the following lemma, where
\begin{equation}\label{eq4}
v^*(x) = [\nabla_{y}^2 g(x, y^*(x))]^{-1} \nabla_{y}f(x, y^*(x)).
\end{equation}

\begin{lemma} \label{lemma:2}
Suppose Assumptions \ref{assum:1} and \ref{assum:2} hold. Then, $\nabla \Phi(x)$ in (\ref{eq2}), $y^*(x)$ in problem (\ref{eq1}), and $v^*(x)$ in  (\ref{eq4}) are respectively $L_{\Phi}$-, $C_{y}$-, and $C_{v}$-Lipschitz continuous, where 
\begin{align*}
&L_{\Phi} = L_f + \frac{2L_fL_{gg}+ C_f^2L_{gg}}{\mu}+ \frac{L_fL_g^2 + 2C_fL_gL_{gg}}{\mu^2} + \frac{C_f L_g^2 L_{gg}}{\mu^3},\\
& C_{y} = \frac{L_g}{\mu},~~~~ C_{v} = \left(\frac{L_f}{\mu} + \frac{C_f L_{gg}}{\mu^2}\right)\left(1 + \frac{L_g}{\mu}\right),~~~~ L_{gg} = \max\{L_{gxy}, L_{gyy}\}.
\end{align*}
\end{lemma}

The proof of Lemma \ref{lemma:2} can be found in \cite[Lemma 2.2]{16} and \cite[Lemma A.1]{56}. To achieve a further improvement in complexity, the following assumption is imposed on the stochastic component functions in problem (\ref{eq1}), as also aopted in \cite{21}.

\begin{assumption}\label{assum:4}
For any $\xi$ and $\zeta$, functions $F(x, y; \xi)$ and $G(x, y; \zeta)$ satisfy the same conditions as in Assumption \ref{assum:2}.
\end{assumption}

In addition, considering the differentiability and nonconvexity of $\Phi(x)$ derived from Assumption \ref{assum:1}, in line with prior works such as \cite{7, 13, 16, 21}, we adopt the $\epsilon$-stationary point defined below as the convergence criterion for problem (\ref{eq1}).

\begin{definition} \label{def:2}
A point $\bar{x}$ is said to be an $\epsilon$-stationary point of problem (\ref{eq1}) if $\mathbb{E}\| \nabla \Phi(\bar{x}) \|^2 \le \epsilon$.
\end{definition}

\section{A SPIDER-based alternating variance-reduced algorithm} \label{sec:3}

In this section, we propose a new alternating stochastic variance-reduced algorithm with multi-step LL updates, based on the SPIDER variance reduction technique  (ALS-SPIDER), and provide its convergence and complexity analysis.

As mentioned in Section \ref{sec:1}, the complexity of the SPIDER-based variance-reduced algorithm VRBO \cite{21}, which employs multi-step LL udpates, does not match that of SPIDER for nonconvex single-level problems. The main reason is the additional bias introduced when approximating the inverse of a Hessian matrix in the estimation of $\nabla \Phi(x)$. Recent studies show that AmIGO \cite{32} and SOBA \cite{34} achieve a compleixty comparable to SGD for nonconvex single-level problems, using a constant number of SGD steps to approximately compute a Hessian-inverse-vector product for estimating $\nabla \Phi(x)$. Additionally, BSVRB$^{\text{v}2}$ \cite{61}, with a one-step LL update, achieves the same complexity as STORM for nonconvex single-level problems by performing a one-step STORM-based update that approximates the Hessian-inverse-vector product. Inspired by recent progress, our main idea is to apply multi-step updates on an auxiliary variable to estimate $\nabla \Phi(x)$. While a one-step update is also feasible, we choose multi-step updates to make the algorithm more flexible and applicable.

Before presenting our algorithm, we first introduce an auxiliary variable $v$ and three directions: $D_{y}(x, y)$, $D_{v}(x, y, v)$, and $D_{x}(x, y, v)$ \cite{32, 34, 35, 22, 43, 61, 59}, which will be useful for the subsequent algorithmic design. These directions are defined as
\begin{align}
&~D_{y}(x, y) = \nabla_{y} g(x, y), ~~~~ D_{v}(x, y, v) = \nabla_{y}^2 g(x, y) v - \nabla_{y} f(x, y), \label{eq5} \\
&D_{x}(x, y, v) = \nabla_{x} f(x, y) - \nabla_{x}\nabla_{y} g(x, y) v.\label{eq6}
\end{align}
Applying gradient descent steps to the LL variable and moving in the direction of $-D_{y}(x, y)$ allows us to approximate $y^*(x)$ in (\ref{eq2}). Once an approximation for $y^*(x)$ is available, $v^*(x)$ in (\ref{eq4}) can be approximated by solving the quadratic programming problem
\begin{equation}\label{eq7}
\min\limits_{v\in \mathbb{R}^q}  \frac{1}{2} v^\top \nabla_{y}^2 g(x, y) v - v^\top \nabla_{y}f(x, y)
\end{equation}
with $v$ following the direction $- D_{v}(x, y, v)$. Notice that when $y^*(x)$ is known exactly, the objective function of problem (\ref{eq7}) takes the form $\frac{1}{2} v^\top \nabla_{y}^2 g(x, y^*(x)) v - v^\top \nabla_{y}f(x, y^*(x))$, and $v^*(x)$ represents the optimal solution. Replacing $y^*(x)$ and $[\nabla_{y}^2 g(x, y^*(x))]^{-1} \nabla_{y}f(x, y^*(x))$ (= $v^*(x)$) in (\ref{eq2}) with the obtained approximate values, $D_{x}(x, y, v)$ serves as a surrogate for $\nabla \Phi(x)$ in (\ref{eq2}) to update $x$.

In the rest of this paper, we focus on the stochastic algorithms for problem (\ref{eq1}). To do so, let 
\begin{equation}\label{eq8}
\bar{\xi}_i^{x} := \{\xi_i^{x}, \delta_i^{x}\}, \qquad \bar{\xi}_i^{y} := \{\delta_i^{y}\}, \qquad \bar{\xi}_i^{v} := \{\xi_i^{v}, \delta_i^{v}\}
\end{equation}
where $\xi_i^{x}$, $\xi_i^{v}$ are samples drawn from the distribution of $\xi$, and $\delta_i^{x}$, $\delta_i^{y}$, $\delta_i^{v}$ are samples drawn from the distribution of $\delta$. We consider the stochastic versions of $D_{y}(x, y)$, $D_{v}(x, y, v)$, and $D_{x}(x, y, v)$ given by 
\begin{align}
&D_{y}(x, y; \bar{\xi}_i^{y}) = \nabla_{y} G(x, y; \delta_i^{y}),  \label{eq9}\\
&D_{v}(x, y, v; \bar{\xi}_i^{v}) = \nabla_{y}^2 G(x, y; \delta_i^{v}) v - \nabla_{y} F(x, y; \xi_i^{v}), \label{eq10}\\
&D_{x}(x, y, v; \bar{\xi}_i^{x}) = \nabla_{x} F(x, y; \xi_i^{x}) - \nabla_{x}\nabla_{y} G(x, y; \delta_i^{x}) v. \label{eq11}
\end{align}
Moreover, we introduce a set $\Omega$ defined as follows.

\begin{definition}\label{definition:3}
We define $\Omega$ as the set $\Omega := \{v \in \mathbb{R}^q: \|v\| \le M\}$, where the radius $M$ satisfies $M \ge C_f/\mu$, and $\mu$, $C_f$ are given in Assumptions \ref{assum:1}, \ref{assum:2}.
\end{definition}

\begin{remark} \label{remark:1} 
Under Assumptions \ref{assum:1}, \ref{assum:2}, and Definition \ref{definition:3}, it is easily obtained that for any $x$, $v^*(x)$ in (\ref{eq4}) belongs to the set $\Omega$.
This result enables us to leverage the non-expansive property of projection to analyze the iterations for updating $v$, as shown in Lemma \ref{lemma:4}.
\end{remark}

\subsection{The proposed algorithm ALS-SPIDER}

ALS-SPIDER alternately updates the UL variable $x$ and the LL variable $y$, and the auxiliary variable $v$ is introduced for estimating $\nabla \Phi(x)$. Multi-step updates are performed on $y$ and $v$ to approximate $y^*(x)$ in (\ref{eq2}) and $v^*(x)$ in (\ref{eq4}), respectively, before each update of $x$. Moreover, the SPIDER variance reduction technique is utilized to construct variance-reduced estimates for updating these variables. Specifically, given $x_{k+1}$, an approximate solution $y_{k+1}$ to the LL problem
\begin{equation}\label{eq13}
\min\limits_{y \in \mathbb{R}^q}g(x_{k+1}, y)
\end{equation}
is obtained by applying Algorithm \ref{alg:2}, which executes $T$ variance-reduced iterations on $y$. Based on the output $y_{k+1}$, Algorithm \ref{alg:3} runs $J$ projected variance-reduced itrations on $v$ to solve the quadratic programming problem (QP)
\begin{equation}\label{eq14}
\min\limits_{v\in \Omega}  \frac{1}{2} v^\top \nabla_{y}^2 g(x_{k+1}, y_{k+1}) v - v^\top \nabla_{y}f(x_{k+1}, y_{k+1})
\end{equation}
where $\Omega$ is defined in Definition \ref{definition:3}, and returns $v_{k+1}$ as the approximate solution. Notice that here we adopt a different approach from (\ref{eq7}) to constrain the domain of $v$ to be $\Omega$, with the goal of ensuring that the iterative trajectory of $v$ remains within a bounded region when using the variance-reduced estimates to update $v$. Then, a periodic variance-reduced estimate $D_{k+1}^x$ is used to update $x_{k+1}$. ALS-SPIDER is summarized in Algorithm \ref{alg:1}.

\textit{Algorithms \ref{alg:2} and \ref{alg:3}.} For each outer loop iteration $k$, Algorithm \ref{alg:2} performs variance-reduced updates on $y$ to approximately solve problem (\ref{eq13}), and $y_{k+1}$ and $D_{k, T-1}^y$ are returned. Specifically, at each $t\in \{0, \ldots, T-1\}$, Algorithm \ref{alg:2} constructs the SPIDER-based variance-reduced estimate 
\begin{equation}\label{eq15}
{D}_{k, t}^{y}= {D}_{y}(\bar{x}_{k, t}, \bar{y}_{k, t}; \mathcal{S}_{2, y}) + (1 - \tau^{y})\left({D}_{k, t-1}^{y}- {D}_{y}(\bar{x}_{k, t-1}, \bar{y}_{k, t-1}; \mathcal{S}_{2, {y}})\right),
\end{equation}
based on the current gradient estimate ${D}_{y}(\bar{x}_{k, t}, \bar{y}_{k, t}; \mathcal{S}_{2, y})$ and the previous gradient estimate ${D}_{y}(\bar{x}_{k, t-1}, \bar{y}_{k, t-1}; \mathcal{S}_{2, {y}})$, using the sample set $\mathcal{S}_{2, y}$ (in line \ref{line:2-3}). Here, $\tau^y$ is set to $0$, and for each $p\in \{t-1, t\}$, ${D}_{y}(\bar{x}_{k, p}, \bar{y}_{k, p}; \mathcal{S}_{2, y})$ takes the same form as in (\ref{eq16}), with $x_{k}$, $y_{k}$, and $S_1$ replaced by $\bar{x}_{k, p}$, $\bar{y}_{k, p}$, and $S_2$, respectively. The variance-reduced estimate $D_{k, t}^y$ is then used to update $\bar{y}_{k, t}$ in line \ref{line:2-2} of Algorithm \ref{alg:2}, while $\bar{x}_{k, t}$ remains unchanged and is set to $x_{k+1}$ (see line \ref{line:9} in Algorithm \ref{alg:1} and line \ref{line:2-4} in Algorithm \ref{alg:2}). At the end, it sets $y_{k+1}$ to $\bar{y}_{k, T-1}$, and returns both $y_{k+1}$ and the recursive variance reduced estimate $D_{k, T-1}^{y}$ as output. Notice that $\bar{x}_{k, t-1}$ and $\bar{y}_{k, t-1}$ at $t=0$ (i.e., $\bar{x}_{k, -1}$ and $\bar{y}_{k, -1}$) are respectively set to $x_k$ and $y_k$, which are obtained from the previous outer loop iteration $k-1$; $D_{k, t-1}^{y}$ at $t=0$ (i.e., $D_{k, -1}^y$) is set to the periodic variance-reduced estimate $D_k^y$ (see line \ref{line:9} in Algorithm \ref{alg:1}). Here, $D_k^y$ is defined as 
\begin{equation}\label{eq16}
D_k^y := D_{y}(x_{k}, y_{k}; \mathcal{S}_{1,y}):= \frac{1}{S_1}\sum\nolimits_{i=1}^{S_1} D_{y}(x_{k}, y_{k}; \bar{\xi}_i^{y})
\end{equation}
if $\mod(k, q_1)=0$, where the gradient estimate $D_{y}(x_{k}, y_{k}; \mathcal{S}_{1,y})$ is computed over the sample set $\mathcal{S}_{1,y}$ (in line \ref{line:3} of Algorithm \ref{alg:1}), $\bar{\xi}_i^{y}$ has the same definition as in (\ref{eq8}), and the computation of $D_{y}(x_k, y_k; \bar{\xi}_i^{y})$ follows (\ref{eq9}).
Otherwise, it sets $D_{k}^{y}$ to $D_{k-1, T-1}^y$, the output of Algorithm \ref{alg:2} in the previous outer loop iteration $k-1$. Such a setting for $\bar{x}_{k, -1}$, $\bar{y}_{k, -1}$, and $D_{k, -1}^{y}$ allows us to track the gradient estimation error $\|D_{k, t}^y - D_{y}(\bar{x}_{k, t}, \bar{y}_{k, t})\|$ to the previous outer loop iteration $k-1$ when $\mod(k, q_1) \neq 0$, thereby establishing a variance reduction property for $\|D_k^y - D_{y}(x_{k}, y_{k})\|$ similar to that of SPIDER for single-level optimization problems.

With the obtained approximation $y_{k+1}$, Algorithm \ref{alg:3} performs variance-reduced updates on $v$ to solve problem (\ref{eq14}), and returns $v_{k+1}$, $D_{k, J-1}^x$, and $D_{k, J-1}^{v}$ as output. Specifically, setting hyperparameters $\tau^x$, $\tau^v$ to $0$. For each $j\in \{0, \ldots, J-1\}$, Algorithm \ref{alg:3} draws two mutually independent sample sets $\mathcal{S}_{2, {x}}$, $\mathcal{S}_{2,{v}}$ (in line \ref{line: 3-3}), and constructs the SPIDER-based variance-reduced estimates $D_{k, j}^{x}$, $D_{k, j}^{v}$, given by  
\begin{align}\label{eq17}
{D}_{k, j}^{x}= {D}_{x}(x_{k, j}, y_{k, j}, v_{k, j}; \mathcal{S}_{2, x}) + (1 - \tau^{x})\left({D}_{k, j-1}^{x}- {D}_{x}(x_{k, j-1}, y_{k, j-1}, v_{k, j-1}; \mathcal{S}_{2, x})\right)\\ \nonumber
{D}_{k, j}^{v}= {D}_{v}(x_{k, j}, y_{k, j}, v_{k, j}; \mathcal{S}_{2, v}) + (1 - \tau^{v})\left({D}_{k, j-1}^{v}- {D}_{v}(x_{k, j-1}, y_{k, j-1}, v_{k, j-1}; \mathcal{S}_{2, v})\right)
\end{align} 
in a way similar to $D_{k, t}^{y}$ in (\ref{eq15}). Here, for each $p\in \{j-1, j\}$, ${D}_{x}(x_{k, p}, y_{k, p}, v_{k, p}; \mathcal{S}_{2, x})$ and ${D}_{v}(x_{k, p}, y_{k, p}, v_{k, p}; \mathcal{S}_{2, v})$ take the same forms as in (\ref{eq20}) and (\ref{eq21}), respectively, with $x_{k}$, $y_{k}$, $v_k$, and $S_1$ replaced by $x_{k, p}$, $y_{k, p}$, $v_{k, p}$, and $S_2$. Then, a projected variance-reduced iteration is applied to update $v_{k, j}$ using $D_{k, j}^{v}$ in line \ref{line: 3-2} of Algorithm \ref{alg:3}, where the projection operator $\text{proj}_{\Omega}(\cdot)$ takes the form of 
\begin{equation}
\text{proj}_{\Omega}(v) =
\begin{cases}
v &\quad\text{if}~~\|v \| \le M, \\
M \frac{v}{\|v\|} &\quad\text{if} ~~\|v\| > M \\
\end{cases}
\label{eq18}
\end{equation}
based on the definition of $\Omega$ in Definition \ref{definition:3}. In constrast, $x_{k, j}$ and $y_{k, j}$ remain unchanged and are respectively set to $x_{k+1}$ and $y_{k+1}$ (see line \ref{line:10} in Algorithm \ref{alg:1} and line \ref{line: 3-5} in Algorithm \ref{alg:3}). After $J$ iterations, Algorithm \ref{alg:3} sets $v_{k+1}$ to $v_{k, J-1}$, and outputs  $v_{k+1}$ and the recursive variance-reduced estimates $D_{k, J-1}^x$, $D_{k, J-1}^{v}$. Notice that, similar to Algorithm \ref{alg:2}, $x_{k, j-1}$, $y_{k, j-1}$, and $v_{k, j-1}$ at $j=0$ (i.e., $x_{k, -1}$, $y_{k, -1}$, and $v_{k, -1}$) are respectively set to $x_k$, $y_k$, and $v_k$ from the previous outer loop iteration $k-1$. In addition, it sets $D_{k, j-1}^{x}$ and $D_{k, j-1}^{v}$ at $j=0$ (i.e., $D_{k, -1}^{x}$ and $D_{k, -1}^{v}$) to the periodic variance-reduced estimates $D_k^x$ and $D_k^v$, respectively. Here, if $\mod(k, q_1)=0$, one has  
\begin{align} 
&D_k^x := D_{x}(x_{k}, y_{k}, v_{k}; \mathcal{S}_{1,x}):= \frac{1}{S_1}\sum\nolimits_{i=1}^{S_1} D_{x}(x_{k}, y_{k}, v_{k}; \bar{\xi}_i^{x}), \label{eq20}\\ 
&D_k^v := D_{v}(x_{k}, y_{k}, v_{k}; \mathcal{S}_{1,v}):= \frac{1}{S_1}\sum\nolimits_{i=1}^{S_1} D_{v}(x_{k}, y_{k}, v_{k}; \bar{\xi}_i^{v}) \label{eq21}
\end{align}
where $\mathcal{S}_{1,x}$ and $\mathcal{S}_{1, v}$ are sampled sets in line \ref{line:3} of Algorithm \ref{alg:1}, and $D_{x}(x_{k}, y_{k}, v_{k}; \bar{\xi}_i^{x})$ and $D_{v}(x_{k}, y_{k}, v_{k}; \bar{\xi}_i^{v})$ follow the definitions in (\ref{eq11}) and (\ref{eq10}), respectively. Otherwise, $D_k^x$ and $D_k^y$ are set to the outputs $D_{k-1, J-1}^x$ and $D_{k-1, J-1}^v$ of Algorithm \ref{alg:3} in the previous outer loop iteration $k-1$, respectively. Such a setting for $x_{k, -1}$, $y_{k, -1}$, $v_{k, -1}$, $D_{k, -1}^x$, and $D_{k,-1}^v$ enables us to establish the variance reduction properties for $\|D_k^x - D_x(x_k,y_k, v_k)\|$ and $\|D_k^v - D_v(x_k,y_k, v_k)\|$, similar to that of SPIDER for single-level optimization problems.

\begin{algorithm}[htb]   
\caption{ALS-SPIDER algorithm for problem (\ref{eq1})} \label{alg:1}
\begin{algorithmic}[1]
\REQUIRE
 $K$, $T$, $J$, $q_1$, $S_1$, $S_2$, stepsizes $\lambda_1$, $\lambda_2$, $\alpha$, $\beta$, $\eta$, initializers ${x}_0\in \mathbb{R}^p$, ${y}_0 \in \mathbb{R}^q$, ${v}_0\in \Omega$; set $\tau^{x}= 0$, $\tau^{y} = 0$, $\tau^{v} = 0$
\FOR{$k = 0,\ldots, K-1$}
\IF{$\mod(k, q_1)=0$}
\STATE Draw $\mathcal{S}_{1,x}:= \{\bar{\xi}_i^{x}, i=1, \ldots, S_1\}$, $\mathcal{S}_{1, y}:= \{\bar{\xi}_i^{y}, i=1, \ldots, S_1\}$, $\mathcal{S}_{1, v}:= \{\bar{\xi}_i^{v}, i=1, \ldots, S_1\}$ with mutually independent samples  \label{line:3}
\STATE Set ${D}_k^{x}$ by (\ref{eq20}), ${D}_k^{y}$ by (\ref{eq16}), and ${D}_k^{v}$ by (\ref{eq21}) \label{line:4}
\ELSE
\STATE Set ${D}_k^{x} = {D}_{k-1, J-1}^{x}$, ${D}_k^{y} = {D}_{k-1, T-1}^{y}$, ${D}_k^{v} = {D}_{k-1, J-1}^{v}$   \label{line:5}
\ENDIF
\STATE Update ${x}_{k+1} = {x}_k - \alpha {D}_{k}^{x}$  \label{line:8}
\STATE Calling Algorithm \ref{alg:2} with $(\bar{x}_{k, -1}, \bar{y}_{k, -1}, {D}_{k, -1}^{y}, \bar{x}_{k, 0})\leftarrow ({x}_{k}, {y}_{k}, {D}_{k}^{y}, {x}_{k+1})$, $(S_2, T, \lambda_1, \beta, \tau^{y})\leftarrow (S_2, T, \lambda_1, \beta, \tau^{y})$, and obtain $({y}_{k+1}, D_{k, T-1}^{y})$ \label{line:9}
\STATE Calling Algorithm \ref{alg:3} with $({x}_{k, -1}, {y}_{k, -1}, {v}_{k, -1}, {D}_{k, -1}^{x}, {D}_{k, -1}^{v}, {x}_{k, 0}, {y}_{k, 0})\leftarrow ({x}_{k}, {y}_{k}, {v}_{k}, {D}_{k}^{x}, {D}_{k}^{v}, {x}_{k +1}, {y}_{k +1})$, $(S_2, J, \lambda_2, \eta, \tau^{x}, \tau^{v})\leftarrow (S_2, J, \lambda_2, \eta, \tau^{x}, \tau^{v})$, and obtain $(v_{k+1}, D_{k, J-1}^{x}, D_{k, J-1}^{v})$\label{line:10}
\ENDFOR
\end{algorithmic}
\end{algorithm}

\begin{remark} 
There are some key differences between ALS-SPIDER and VRBO \cite{21} in terms of algorithmic design. On one hand, VRBO requires a three-level loop structure, whereas, as shown in Algorithm \ref{alg:1}, ALS-SPIDER features a two-level loop structure. On the other hand, VRBO uses Neumann series to estimate $\nabla \Phi(x)$, while ALS-SPIDER estimates $\nabla \Phi(x)$ via variance-reduced iterations in Algorithm \ref{alg:3}. This difference complicates the analysis of ALS-SPIDER, as it involves handling additional recursive solution errors and variance reduction errors for estimating $\nabla \Phi(x)$.

\end{remark}

\begin{algorithm}[htb]   
\caption{Multi-step variance-reduced iterations for solving the LL problem in (\ref{eq1})} \label{alg:2}
\begin{algorithmic}[1]
\REQUIRE
$\bar{x}_{k, -1}$, $\bar{y}_{k, -1}$, ${D}_{k, -1}^{y}$, $\bar{x}_{k, 0}$, $S_2$, $T$, $\lambda_1$, $\beta$, $\tau^{y}$
\FOR {$t = 0, \ldots, T-1$}
\STATE Update $\bar{y}_{k, t} = \bar{y}_{k, t-1} - \lambda_1 \beta {D}_{k, t-1}^{y}$   \label{line:2-2}
\STATE Draw $\mathcal{S}_{2, y}:= \{\bar{\xi}_i^{y}, i=1, \ldots, S_2\}$ with mutually independent samples \label{line:2-3}
\STATE Compute ${D}_{k, t}^{y}$ by (\ref{eq15}), and then set $\bar{x}_{k, t+1} = \bar{x}_{k, t}$ \label{line:2-4}
\ENDFOR
\STATE Set $y_{k+1}= \bar{y}_{k, T-1}$
\STATE \textbf{return} $(y_{k+1}, {D}_{k, T-1}^{y})$
\end{algorithmic}
\end{algorithm}

\begin{algorithm}[htb]   
\caption{Multi-step projected variance-reduced iterations for appximating $v^*(x)$ in (\ref{eq4})} \label{alg:3}
\begin{algorithmic}[1]
\REQUIRE
${x}_{k, -1}$, ${y}_{k, -1}$, ${v}_{k, -1}$, ${D}_{k, -1}^{x}$, ${D}_{k, -1}^{v}$, ${x}_{k, 0}$, ${y}_{k, 0}$, $S_2$, $J$, $\lambda_2$, $\eta$, $\tau^{x}$, $\tau^{v}$
\FOR {$j = 0, \ldots, J-1$} 
\STATE Update ${v}_{k, j} = \text{proj}_{\Omega}({v}_{k, j-1} -\lambda_2 \eta {D}_{k, j-1}^{v})$, with $\text{proj}_{\Omega}(\cdot)$ defined in (\ref{eq18})   \label{line: 3-2}
\STATE Draw $\mathcal{S}_{2, {x}}:= \{\bar{\xi}_i^{x}, i=1, \ldots, S_2\}$, $\mathcal{S}_{2,{v}}:= \{\bar{\xi}_i^{v}, i=1, \ldots, S_2\}$ with mutually independent samples \label{line: 3-3}
\STATE Compute both $D_{k,j}^x$ and ${D}_{k, j}^{v}$ by (\ref{eq17}) \label{line: 3-4}
\STATE Set ${x}_{k, j+1} = {x}_{k, j}$, ~~${y}_{k, j+1} = {y}_{k, j}$  \label{line: 3-5}
\ENDFOR  
\STATE Set ${v}_{k+1} ={v}_{k, J-1}$   
\STATE \textbf{return} $(v_{k+1}, {D}_{k, J-1}^{x}, {D}_{k, J-1}^{v})$
\end{algorithmic}
\end{algorithm}

\subsection{Theoretical Analysis}

In this section, we begin by presenting several lemmas to describe the properties of the algorithm. Leveraging these lemmas, we then provide the convergence and complexity analysis of ALS-SPIDER.

\subsubsection{Fundamental Lemmas}

To analyze the convergence of ALS-SPIDER, we first examine the change in the function value $\Phi(x)$ after a single update step on the UL variale $x$, as executed in line \ref{line:8} of Algorithm \ref{alg:1}.

\begin{lemma} \label{lemma:3}
Apply Algorithm \ref{alg:1} to solve problem (\ref{eq1}) under Assumptions \ref{assum:1} and \ref{assum:2}. Define $\Delta \Phi_k := \mathbb{E}[\Phi(x_{k+1})-\Phi(x_{k})]$. Then for any non-negative integer $k$,
\begin{align} \label{eq22}
\Delta \Phi_k \le -\frac{\alpha}{2}\mathbb{E}\|\nabla \Phi(x_k)\|^2 - \frac{\alpha}{2}(1 - \alpha L_{\Phi})\mathbb{E}\|D_{k}^{x}\|^2 + 4L_g^2 \alpha \delta_{k}^{v} + \bar{M}_1\alpha \delta_k^{y} + \alpha V_k^{x}
\end{align}
where $\delta_k^{y}= \mathbb{E}\|y_k - y^*(x_k)\|^2$, $\delta_k^{v}= \mathbb{E}\|v_k - v^*(x_k)\|^2$, $V_k^{x}= \mathbb{E}\|D_k^{x}- D_{x}(x_k, y_k, v_k)\|^2$, $\bar{M}_1= 4L_f^2 + 4M^2L_{gg}^2$. Here, $L_{\Phi}$ and $L_{gg}$ are in Lemma \ref{lemma:2}, $M$ is in Definition \ref{definition:3}, and $D_{x}(x_k, y_k, v_k)$ and $v^*(x_k)$ follow the definitions in  (\ref{eq6}) and (\ref{eq4}), respectively.
\end{lemma}
\begin{proof}
Based on the Lipschitz continuity of $\nabla \Phi(x)$ in Lemma \ref{lemma:2}, and step \ref{line:8} of Algorithm \ref{alg:1}, we have, for any non-negative integer $k$, 
\begin{align} \label{eq25}
&\Phi(x_{k+1})- \Phi(x_k) \\ \nonumber
&\le  \langle \nabla \Phi(x_k), x_{k+1} - x_k \rangle + \frac{L_{\Phi}}{2}\| x_{k+1} - x_k \|^2  \\ \nonumber
& = - \frac{\alpha}{2}\|\nabla \Phi(x_k)\|^2 - \frac{\alpha}{2}(1 - \alpha L_{\Phi})\|D_k^{x}\|^2 + \frac{\alpha}{2}\| \nabla \Phi(x_k) - D_k^{x}\|^2 \\ \nonumber
& \le - \frac{\alpha}{2} \|\nabla \Phi(x_k)\|^2 - \frac{\alpha}{2}(1 - \alpha L_{\Phi})\|D_k^{x}\|^2 \\ \nonumber
&~~~ + {\alpha}\|D_k^{x} - D_{x}(x_k, y_k, v_k)\|^2  +  \alpha\|\nabla \Phi(x_k) - D_{x}(x_k, y_k, v_k)\|^2.
\end{align}
Using (\ref{eq2}), (\ref{eq4}), and (\ref{eq6}), we obtain
\begin{align}\label{eq26}
\|\nabla \Phi(x_k) - D_{x}(x_k, y_k, v_k)\|^2 \le C_1+ C_2
\end{align}

where $C_1 = 2\|\nabla_{x}f(x_k, y_k) - \nabla_{x}f(x_k, y^*(x_k))\|^2 + 4\|\nabla_{x}\nabla_{y}g(x_k, y_k)\|^2 \|v_k - v^*(x_k)\|^2$, $C_2 = 4\|v^*(x_k)\|^2 \|\nabla_{x}\nabla_{y}g(x_k, y_k) - \nabla_{x}\nabla_{y}g(x_k, y^*(x_k))\|^2 $. From Assumptions \ref{assum:2}, $L_{gg}$ in Lemma \ref{lemma:2}, and $M$ in Definition \ref{definition:3} ($M \ge C_f/\mu$), it holds that
\begin{align}\label{eq26-1}
C_1 + C_2 \le (2L_f^2 + 4M^2 L_{gg}^2)\|y_k - y^*(x_k)\|^2 + 4L_g^2 \|v_k - v^*(x_k)\|^2.
\end{align} 
Then, the proof is completed by combining (\ref{eq25}), (\ref{eq26}), and (\ref{eq26-1}), utilizing the definition of $\bar{M}_1$, and taking a full expectation.
\end{proof}

Lemma \ref{lemma:3} shows that the upper bound for $\Delta \Phi_k$ depends on the hypergradient estimation variance $V_k^{x}$, the LL problem approximation error $\delta_k^{y}$, and the QP solution error $\delta_k^{v}$. Let $i$, $k$, $T$, $J$ be arbitrary. Using $D_y(x, y)$ and $D_{v}(x, y, v)$ in (\ref{eq5}), we introduce the following notation for Algorithms \ref{alg:1}, \ref{alg:2}, and \ref{alg:3}
\begin{align}\label{eq26-2}
&V_k^y = \mathbb{E}\|D_k^y - D_y(x_k, y_k)\|^2,~ V_k^v = \mathbb{E}\|D_k^v - D_v(x_k, y_k, v_k)\|^2, \\ \nonumber
&W_{i,T}^{y} = \sum\nolimits_{t=0}^{T-1}\mathbb{E}\|D_{i, t-1}^{y}\|^2, ~~V_{i, T}^{y} = \sum\nolimits_{t=0}^{T-1}\mathbb{E}\|D_{i, t-1}^{y} - D_y (\bar{x}_{i, t-1},\bar{y}_{i, t-1})\|^2, \\ \nonumber
&W_{i,J}^{v} = \sum\nolimits_{j=0}^{J-1}\mathbb{E}\|D_{i, j-1}^{v}\|^2, V_{i, J}^{v}= \sum\limits_{j=0}^{J-1}\mathbb{E}\|D_{i, j-1}^{v} - D_{v}(x_{i, j-1}, y_{i, j-1}, v_{i, j-1})\|^2,
\end{align} 
where $V_k^y$ and $V_k^v$ denote the variances at the $k$-th outer iteration, while $W_{i,T}^{y}$, $V_{i,T}^{y}$ (resp. $W_{i,J}^{v}$, $V_{i,J}^{v}$) represent the cumulative gradient estimates and cumulative bias over $T$(resp. $J$) iterations for solving the LL problem(resp. quadratic programming problem). Then, inequalities for $\delta_k^{y}$ and $\delta_k^{v}$ are derived in Lemma \ref{lemma:4}, and for $V_{k}^{x}$, $V_{k, T}^{y}$, and $V_{k, J}^{v}$ are provided in Lemma \ref{lemma:5}. 

\begin{lemma}\label{lemma:4} 
Suppose Assumptions \ref{assum:1} and \ref{assum:2} hold. Define $r_{y}= \mu \lambda_1 \beta/(4 - 2\mu \lambda_1 \beta)$, $r_{v} = \mu \lambda_2 \eta/(4 - 2\mu \lambda_2 \eta)$. Let $0< \beta <1$, $0< \eta <1$, $0< \lambda_1 \le 1/(2L_g + 2\mu)$, $0< \lambda_2 \le 1/(2L_g + 2\mu)$, and $\sum_{i=s}^{s-1}(\cdot)=0$ for any non-negative integer $s$. Then, for any non-negative integer $k$, and positive integers $T$ and $J$, 
\begin{align*}  
&\delta_{k+1}^{y} \le \beta_{yy}\delta_k^{y}+  \beta_{yx}\mathbb{E}\|D_k^{x}\|^2 - (1/2)^{2T-1} \lambda_1^2 \beta W_{k,T}^{y} + \beta_{T}^{y}V_{k, T}^{y}, \\ \nonumber
& \delta_{k+1}^{v} \le \beta_{vv}\delta_k^{v} + \beta_{vx}\mathbb{E}\|D_k^{x}\|^2 - (1/2)^{2J-1}\lambda_2^2 \eta W_{k,J}^{v} + \beta_{J}^{v}V_{k, J}^{v} + C_1(J)(\delta_k^{y}+  \delta_{k+1}^{y})
\end{align*}
where $\beta_{yy}=\left(1 - {\mu \lambda_1\beta}/{2}\right)^T(1 + r_{y})$, $\beta_{T}^{y}= {8\lambda_1\beta}/{\mu}$, $\beta_{yx}= (1 + 1/r_{y})C_{y}^2 \alpha^2 + L_g^2 \alpha^2\beta_{T}^{y}$, $\beta_{vv}=\left(1 - {\mu \lambda_2\eta}/{2}\right)^J(1 + r_{v})$, $\beta_{J}^{v}= {8 \lambda_2 \eta}/{\mu}$, $\beta_{vx} = (1 + 1/r_{v})C_{v}^2 \alpha^2 + \bar{M}_1(1 + 2C_{y}^2)\alpha^2\beta_{J}^{v}$, and $C_1(J)= \bar{M}_1J\beta_{J}^{v}$. Here, $C_{y}$ and $C_{v}$ are in Lemma \ref{lemma:2}, $\bar{M}_1$, $\delta_k^{y}$, and $\delta_k^{v}$ are in Lemma \ref{lemma:3}, and $W_{k,T}^{y}$, $W_{k,J}^{v}$, $V_{k, T}^{y}$, and $V_{k, J}^{v}$ are in (\ref{eq26-2}) with $i=k$.
\end{lemma}

\begin{proof}
We first prove the inequality for $\delta_{k+1}^{y}$. Let the non-negative integer $k$ be arbitrary. Define $\tilde{y}_{k, t} = \bar{y}_{k, t-1} -\lambda_1 D_{k, t-1}^{y}$. From step \ref{line:2-2} of Algorithm \ref{alg:2}, we obtain
\begin{align}\label{eq27}  
\|\bar{y}_{k, t} - y^*(x_{k+1})\|^2 &= \|\bar{y}_{k, t-1} - \lambda_1 \beta D_{k, t-1}^{y} - y^*(x_{k+1})\|^2 \\ \nonumber 
&= \|\bar{y}_{k, t-1} - y^*(x_{k+1})\|^2 + \lambda_1^2 \beta^2 \|D_{k, t-1}^{y}\|^2 \\ \nonumber
&~~~~- 2 \lambda_1^2 \beta \|D_{k, t-1}^{y}\|^2 - 2 \lambda_1 \beta \langle D_{k, t-1}^{y}, \tilde{y}_{k, t} - y^*(x_{k+1})\rangle.
\end{align}

By the $\mu$-strong convexity and $L_g$-smoothness of $g(x_{k+1}, \cdot)$ w.r.t.$y$, for any $y \in \mathbb{R}^q$, 
\begin{align}\label{eq28}
&\langle D_{k, t-1}^{y}, y - \tilde{y}_{k, t}\rangle \\ \nonumber
&\le g(x_{k+1}, y) - g(x_{k+1}, \bar{y}_{k, t-1})  - \langle \nabla_{y}g(x_{k+1}, \bar{y}_{k, t-1}) - D_{k, t-1}^{y}, y- \tilde{y}_{k, t} \rangle \\ \nonumber
&~~~~ - \langle \nabla_{y}g(x_{k+1}, \bar{y}_{k, t-1}),\tilde{y}_{k, t}- \bar{y}_{k, t-1} \rangle - \frac{\mu}{2}\| y - \bar{y}_{k, t-1} \|^2\\ \nonumber
& \le g(x_{k+1}, y)-  g(x_{k+1}, \tilde{y}_{k, t})- \langle \nabla_{y}g(x_{k+1}, \bar{y}_{k, t-1}) - D_{k, t-1}^{y}, y - \tilde{y}_{k, t} \rangle\\ \nonumber
&~~~~- \frac{\mu}{2}\|y - \bar{y}_{k, t-1} \|^2 + \frac{L_g}{2}\|\tilde{y}_{k, t}- \bar{y}_{k,t-1}\|^2.
\end{align}

Applying Young's inequality, it holds that
\begin{align}\label{eq30}
&\langle \nabla_{y}g(x_{k+1}, \bar{y}_{k, t-1}) - D_{k, t-1}^{y}, y - \tilde{y}_{k, t} \rangle\\ \nonumber
&\ge -\frac{2}{\mu} \|\nabla_{y}g(x_{k+1}, \bar{y}_{k, t-1}) - D_{k, t-1}^{y}\|^2 - \frac{\mu}{4}(\|y -\bar{y}_{k, t-1}\|^2 + \|\tilde{y}_{k, t} - \bar{y}_{k, t-1}\|^2).
\end{align}

Combining (\ref{eq28}) and  (\ref{eq30}), and letting $y = y^*(x_{k+1})$, one has
\begin{align}\label{eq31}
&\langle D_{k, t-1}^{y}, y^*(x_{k+1})- \tilde{y}_{k, t}\rangle \\ \nonumber
&\le g(x_{k+1}, y^*(x_{k+1}))- g(x_{k+1}, \tilde{y}_{k, t})  + \frac{2}{\mu}\|\nabla_{y}g(x_{k+1}, \bar{y}_{k, t-1}) - D_{k, t-1}^{y} \|^2  \\ \nonumber
& ~~~~- \frac{\mu}{4}\| y^*(x_{k+1}) - \bar{y}_{k, t-1} \|^2 + \big(\frac{L_g}{2} + \frac{\mu}4\big)\|\tilde{y}_{k, t} - \bar{y}_{k, t-1}\|^2 \\ \nonumber
& \le \frac{2}{\mu}\|\nabla_{y}g(x_{k+1}, \bar{y}_{k, t-1}) - D_{k, t-1}^{y} \|^2 - \frac{\mu}{4}\|y^*(x_{k+1}) - \bar{y}_{k, t-1} \|^2 + \frac{\lambda_1}{4} \|D_{k,t-1}^{y}\|^2.
\end{align}
where the last inequality is due to the definition of $\tilde{y}_{k, t}$, the fact that $y^*(x_{k+1})$ minimizes (\ref{eq13}), and the choice of $\lambda_1$. 

It then follows from (\ref{eq27}) and (\ref{eq31}) that, for any non-negative integer $t$, 
\begin{align} \label{eq32} 
\|\bar{y}_{k, t} - y^*(x_{k+1})\|^2 &\le \big(1 - \frac{\mu}{2} \lambda_1\beta\big) \|\bar{y}_{k, t-1} - y^*(x_{k+1})\|^2 - \frac{1}{2}\lambda_1^2 \beta \|D_{k, t-1}^{y}\|^2 \\ \nonumber
&~~~~ + \frac{4}{\mu} \lambda_1 \beta \|\nabla_{y}g(x_{k+1}, \bar{y}_{k, t-1}) - D_{k, t-1}^{y} \|^2.
\end{align}

Since $\|y_{k+1}- y^*(x_{k+1}) \|^2 = \|\bar{y}_{k, T-1} - y^*(x_{k+1})\|^2$, using induction on (\ref{eq32}), and leveraging $\bar{y}_{k, -1}= y_k$, $D_{k, -1}^{y}= D_k^{y}$, $\sum_{t=0}^{-1}(\cdot)=0$, and $\bar{x}_{k, t} = x_{k+1}$ for any non-negative integer $t$, we have, for any positive integer $T$,
\begin{align*} 
\|y_{k+1}- y^*(x_{k+1}) \|^2  & \le (1 - \frac{\mu}{2} \lambda_1 \beta)^{T} \|y_{k} - y^*(x_{k+1})\|^2 +  \frac{4}{\mu}\lambda_1 \beta T_{1, y}\\ \nonumber
&~~+ \frac{4}{\mu}\lambda_1 \beta\sum\limits_{t=0}^{T-2}\|\nabla_{y}g(\bar{x}_{k,t}, \bar{y}_{k,t}) - D_{k,t}^{y} \|^2 - \frac{1}{2}\lambda_1^2 \beta T_{2, y}
\end{align*}
with $T_{1, y} = \|\nabla_{y}g(x_{k+1}, y_k)- D_k^{y}\|^2$, $T_{2, y}= \sum\nolimits_{t=0}^{T-1}(1 -\frac{\mu}{2}\lambda_1 \beta)^{T-1-t}\|D_{k, t-1}^{y}\|^2$.

By Assumption \ref{assum:2}, the Lipschitz continuity of $y^*(x)$ in Lemma \ref{lemma:2}, Young's inequality, step \ref{line:8} of Algorithm \ref{alg:1}, and the choices of $\lambda_1$ and $\beta$, one can see that
\begin{align*}
& \|y_k - y^*(x_{k+1})\|^2 \le (1 + r_{y})\|y_k - y^*(x_k)\|^2 + (1 + 1/r_{y})C_{y}^2 \alpha^2 \|D_k^{x}\|^2,\\ \nonumber
& T_{1, y} \le 2(\|\nabla_{y}g(x_{k}, y_k) - D_k^{y}\|^2 + L_g^2\alpha^2 \|D_k^{x}\|^2), T_{2, y}\ge ({1}/{4})^{T-1}\sum\nolimits_{t=0}^{T-1}\|D_{k, t-1}^{y}\|^2.
\end{align*}
Substituting these relations into the preceding inequality, using (\ref{eq5}), $D_{k, -1}^{y} = D_k^{y}$, $x_k = \bar{x}_{k, -1}$, $y_k = \bar{y}_{k, -1}$, and taking a full expectation, the inequality for $\delta_{k+1}^{y}$ is obtained.

Then, we prove the inequality for $\delta_{k+1}^{v}$. Let the non-negative integer $k$ be arbitrary. Define $\tilde{v}_{k, j} = v_{k, j-1} -\lambda_2 D_{k, j-1}^{v}$. From step \ref{line: 3-2} of Algorithm \ref{alg:3}, the fact that $v^*(x_{k+1}) \in \Omega$, and the non-expansive property of projection,
\begin{align}\label{eq34} 
\|v_{k, j} - v^*(x_{k+1})\|^2 & \le \|v_{k, j-1} - \lambda_2 \eta D_{k, j-1}^{v} - v^*(x_{k+1})\|^2 \\ \nonumber
& = \|v_{k, j-1} - v^*(x_{k+1})\|^2 + \lambda_2^2 \eta(\eta -2)\| D_{k, j-1}^{v}\|^2  \\ \nonumber
& ~~~~- 2 \lambda_2 \eta \langle D_{k, j-1}^{v},   \tilde{v}_{k, j} - v^*(x_{k+1})\rangle.
\end{align}

Let $\bar{G}(v) = v^{\top} \nabla_{y}^2 g(x_{k+1}, y^*(x_{k+1}))v - v^{\top} \nabla_{y}f(x_{k+1}, y^*(x_{k+1}))$. It is easy to verify that $\bar{G}(v)$ is $L_g$-smooth and $\mu$-strongly convex w.r.t.$v$, and $v^*(x_{k+1})$ minimizes $\min_{v \in \mathbb{R}^q}\bar{G}(v)$. With a proof similar to (\ref{eq28}), we get, for any $v \in \mathbb{R}^q$, 
\begin{align}\label{eq35} 
\langle D_{k, j-1}^{v},  v - \tilde{v}_{k, j}\rangle & \le  \bar{G}(v) - \bar{G}( \tilde{v}_{k,j})- \langle \nabla_{v} \bar{G}(v_{k, j-1}) - D_{k, j-1}^{v}, v - \tilde{v}_{k,j}\rangle \\ \nonumber
& ~~~~-  \frac{\mu}{2}\|v - v_{k, j-1}\|^2  + \frac{L_g}{2}\| \tilde{v}_{k,j}- v_{k, j-1}\|^2.
\end{align}

Using Young's inequality, one has
\begin{align}\label{eq35-1} 
&\langle \nabla_{v} \bar{G}(v_{k, j-1}) -  D_{k, j-1}^{v}, v - \tilde{v}_{k,j}\rangle \\ \nonumber
& \ge - \frac{2}{\mu}\|\nabla_{v} \bar{G}(v_{k, j-1}) - D_{k, j-1}^{v}\|^2 - \frac{\mu}{4}(\| v - v_{k,j-1}\|^2 +\|v_{k,j-1}- \tilde{v}_{k,j}\|^2 ).
\end{align}

Combining (\ref{eq35}) and (\ref{eq35-1}), and letting $v= v^*(x_{k+1})$, we obtain 
\begin{align}\label{eq36}
\langle D_{k, j-1}^{v},  v^*(x_{k+1})- \tilde{v}_{k, j}\rangle &\le {2}{\mu}^{-1}\|\nabla_{v} \bar{G}(v_{k, j-1})- D_{k, j-1}^{v}\|^2  \\ \nonumber
&~~~~ - {\mu}/{4}\| v^*(x_{k+1})- v_{k, j-1} \|^2  + {\lambda_2}/{4}\|D_{k, j-1}^{v}\|^2 
\end{align}
following a derivation similar to (\ref{eq31}).

For the first term on the right-hand side (r.h.s.) of (\ref{eq36}), it holds that
\begin{align} \label{eq37}
\|D_{k, j-1}^{v}- \nabla_{v}\bar{G}(v_{k, j-1})\|^2 \le 2\|D_{k, j-1}^{v}- D_{v}(x_{k,j-1}, y_{k, j-1}, v_{k, j-1})\|^2 + T_0.
\end{align}
where $T_0 = 2\|\nabla_{v}\bar{G}(v_{k, j-1}) - D_{v}(x_{k,j-1}, y_{k, j-1}, v_{k, j-1})\|^2$. Applying the definition of $\bar{G}(v)$, (\ref{eq5}), $\|v_{k, j-1}\|\le M$, Assumption \ref{assum:2}, and $\bar{M}_1$ in Lemma \ref{lemma:3}, one can see that
\begin{align} \label{eq38}
T_0 &\le 4M^2 \|\nabla_{y}^2 g(x_{k+1}, y^*(x_{k+1}))- \nabla_{y}^2 g(x_{k, j-1}, y_{k, j-1})\|^2 \\ \nonumber
&~~~+ 4  \|\nabla_{y} f(x_{k+1}, y^*(x_{k+1}))- \nabla_{y} f(x_{k, j-1}, y_{k, j-1})\|^2  \\ \nonumber
& \le \bar{M}_1(\|x_{k+1} - x_{k, j-1}\|^2 + \|y^*(x_{k+1}) - y_{k, j-1}\|^2).
\end{align}

By (\ref{eq34}), (\ref{eq36}), (\ref{eq37}), and (\ref{eq38}), we have, for any non-negative integer $j$,
\begin{align}\label{eq39}
&\|v_{k, j} - v^*(x_{k+1})\|^2 \\ \nonumber
& \le \big(1 - \mu \lambda_2 \eta/2 \big)\|v_{k, j-1} - v^*(x_{k+1})\|^2 - \lambda_2^2 \eta/ {2} \|D_{k, j-1}^{v}\|^2 \\ \nonumber
&~~~+ 4 \lambda_2 \eta \bar{M}_1\mu^{-1}(\|x_{k+1} - x_{k, j-1}\|^2 + \|y^*(x_{k+1}) - y_{k, j-1}\|^2) \\ \nonumber
&~~~+ \beta_J^{v}\|D_{v}(x_{k,j-1}, y_{k, j-1}, v_{k, j-1}) - D_{k, j-1}^{v}\|^2.
\end{align}

Since $\|v_{k+1} - v^*(x_{k+1})\|^2 =  \|v_{k, J-1} - v^*(x_{k+1})\|^2$, using induction on (\ref{eq39}), leveraging $(x_{k, -1}, y_{k,-1}, v_{k, -1})=(x_k, y_k, v_k)$, $x_{k, j}= x_{k+1}$ and $y_{k, j}= y_{k+1}$ for any non-negative integer $j$, and $\sum\nolimits_{j=1}^{0}(\cdot)=0$, it follows that, for any positive integer $J$, 
\begin{align*} 
\|v_{k+1} - v^*(x_{k+1})\|^2 &\le  \big(1 - \mu \lambda_2 \eta/2  \big)^{J}\| v_{k} - v^*(x_{k+1})\|^2 - \lambda_2^2 \eta/2  T_{1, v} \\ \nonumber
&~~~ + {4}\lambda_2\eta\bar{M}_1\mu^{-1} \left(T_{2, v} + (J-1)\|y_{k+1} - y^*(x_{k+1})\|^2\right) \\ \nonumber
& ~~~+ \beta_{J}^{v} \sum\nolimits_{j=0}^{J-1}\|D_{k,j-1}^{v} - D_{v}(x_{k,j-1}, y_{k,j-1}, v_{k,j-1})\|^2
\end{align*}
with $T_{1, v}= \sum\nolimits_{j=0}^{J-1}(1 - \frac{\mu}{2} \lambda_2 \eta )^{J-1-j}\|D_{k, j-1}^{v}\|^2$, $T_{2, v}= \|x_{k+1} - x_k\|^2 + \|y_k - y^*(x_{k+1})\|^2$.

Utilizing Young's inequality, the Lipschitz continuity of $v^*(x)$ and $y^*(x)$ in Lemma \ref{lemma:2}, step \ref{line:8} of Algorithm \ref{alg:1}, and the choices of $\lambda_2$ and $\eta$, it is easy to check that 
\begin{align*}
&\|v_k - v^*(x_{k+1})\|^2 \le (1 + r_{v})\|v_k - v^*(x_k)\|^2 + (1 + 1/r_{v})C_{v}^2\alpha^2 \|D_k^{x}\|^2,\\ \nonumber
& T_{1, v}\ge (1/4)^{J-1}\sum\nolimits_{j=0}^{J-1}\|D_{k, j-1}^{v}\|^2, T_{2, v} \le 2\|y_k - y^*(x_k)\|^2 + (1 + 2C_{y}^2)\alpha^2\|D_k^{x}\|^2.
\end{align*}
Substituting these relations into the preceding inequality, and taking a full expectation, the inequality for $\delta_{k+1}^{v}$ is obtained. 
\end{proof}

\begin{lemma}\label{lemma:5}
Suppose Assumptions \ref{assum:1}, \ref{assum:2}, \ref{assum:3}, and \ref{assum:4} hold. Define $\eta_{y} = {\bar{M}_1}T \lambda_1^2 \beta^2$, $\eta_{v} = 4L_g^2 \lambda_2^2 \eta^2$. Let $0 \le \tau^x < 1$, $0 \le \tau^y < 1$, $0 \le \tau^v < 1$, and $\sum_{i=s}^{s-1}(\cdot) = 0$ for any non-negative integer $s$. Then, for any non-negative integer $k$, positive integers $T$ and $J$, 
\begin{align}
&V_{k}^{x}  \le \frac{2}{S_2}\sum\nolimits_{i=(n_k-1)q_1}^{k-1} W_{i,T,J} + \left(\frac{2}{S_1} + \frac{4}{S_2}\tau^{x}\right) J^2 (\sigma_f^2 + M^2\sigma_{gxy}^2), \label{eq40-1}\\ 
&V_{k, T}^{y} \le \frac{2}{S_2}TL_g^2\sum\nolimits_{i=(n_k-1)q_1}^{k}W_{i, T} + \left(\frac{1}{S_1} +  \frac{4}{S_2}\tau^{y}\right)T^2 \sigma_g^2, \label{eq40-2} \\ 
&V_{k, J}^{v} \le \frac{2}{S_2}J\sum\nolimits_{i=(n_k-1)q_1}^{k}W_{i,T,J} + \left(\frac{2}{S_1}+ \frac{8}{S_2}\tau^v \right) J^2 (\sigma_f^2 + M^2 \sigma_{gyy}^2)\label{eq40-3}
\end{align}
where $W_{i, T}=\alpha^2\mathbb{E}\|D_i^{x}\|^2 + \lambda_1^2 \beta^2 W_{i, T}^{y}$, $W_{i,T,J}= \bar{M}_1 \alpha^2\mathbb{E}\|D_i^{x}\|^2 + \eta_{y} W_{i, T}^{y} + \eta_{v} W_{i, J}^{v}$, and $n_k$ is the integer satisfying $k = (n_k-1)q_1$ when $\mod(k, q_1)=0$, and $(n_k-1)q_1 <k <(n_k) q_1$ when $\mod(k, q_1)\neq 0$. Here, $\bar{M}_1$ and $V_{k}^{x}$ are in Lemma \ref{lemma:3}, $M$ is in Definition \ref{definition:3}, $W_{i,T}^{y}$ and $W_{i,J}^{v}$ are in (\ref{eq26-2}), and $V_{k, T}^{y}$ and $V_{k, J}^{v}$ are in Lemma \ref{lemma:4}.
\end{lemma}

\begin{proof}
For any $k$, $t$, $j$, define $\mathcal{F}_{k, t} = \sigma\{ x_0, y_0, v_0, \ldots, x_k, \ldots, \bar{y}_{k-1, t}\}$, $\mathcal{F}_{k, j} = \sigma\{ x_0, \ldots, x_k, \ldots, y_k, \ldots, v_{k-1, j}\}$, $\mathcal{F}_{k}^{'} = \sigma\{ x_0, \ldots, x_k, \ldots, y_k, \ldots, v_{k-1, J-1}, v_k\}$, where $\sigma\{\cdot\}$ denotes the $\sigma$-algebra generated by random variables. Moreover, we define $u_{k}:=(x_{k}, y_{k}, v_{k})$, $\bar{u}_{k, t}:=(\bar{x}_{k, t}, \bar{y}_{k, t})$, and $u_{k, j}:=(x_{k, j}, y_{k, j}, v_{k, j})$ for simplicity.

Let positive integers $T$ and $J$ be arbitrary. We first prove the inequality for $V_{k}^{x}$. From step \ref{line: 3-4} of Algorithm \ref{alg:3}, we have, for any non-negative integers $k$ and $j$,
\begin{align}\label{eq41}
&\mathbb{E}[\|D_{k, j}^{x} - D_{x}(u_{k, j})\|^2 | \mathcal{F}_{k+1, j}] \\ \nonumber
& = (1 - \tau^x)^2 \|D_{k, j-1}^{x} - D_{x}(u_{k, j-1})\|^2  + \mathbb{E}[P_1|\mathcal{F}_{k+1, j}]
\end{align}
where $P_1 = \|D_{x}(u_{k, j}; \mathcal{S}_{2, x})- D_{x}(u_{k, j}) + (1 - \tau^x) \big(D_{x}(u_{k, j-1})- D_{x}(u_{k, j-1}; \mathcal{S}_{2, x})\big)\|^2$. 

By the triangle inequality, 
\begin{align}\label{eq42-0}
\mathbb{E}[P_1|\mathcal{F}_{k+1, j}] & \le 2(1 - \tau^x)^2 \mathbb{E}[ P_{11} |\mathcal{F}_{k+1, j}] + 2(\tau^x)^2 \mathbb{E}[ P_{12} |\mathcal{F}_{k+1, j}]
\end{align}
where $P_{11} = \|D_{x}(u_{k, j}; \mathcal{S}_{2, x}) - D_{x}(u_{k, j-1}; \mathcal{S}_{2, x}) - D_{x}(u_{k, j}) + D_{x}(u_{k, j-1})\|^2$, $P_{12} = \|D_{x}(u_{k, j}; \mathcal{S}_{2, x}) - D_{x}(u_{k, j}) \|^2$.

Applying the definitions of $D_{x}(u_{k, j}; \mathcal{S}_{2, x})$ and $D_{x}(u_{k, j-1}; \mathcal{S}_{2, x})$, the mutual independence of $\mathcal{S}_{2, x}$, and the inequality $\mathbb{E}[\|x - \mathbb{E}[x| \mathcal{F}_{k+1, j}]\|^2| \mathcal{F}_{k+1, j}]\le \mathbb{E}[\| x \|^2 |\mathcal{F}_{k+1, j}]$, one has
\begin{align}
\mathbb{E}[P_{11}|\mathcal{F}_{k+1, j}] & \le \frac{1}{S_2^2} \sum\nolimits_{i=1}^{S_2}\mathbb{E}[\| D_{x}(u_{k, j}; \bar{\xi}_i^{x}) - D_{x}(u_{k, j-1};\bar{\xi}_i^{x}) \|^2 |\mathcal{F}_{k+1, j}], \label{eq42-1}\\ 
\mathbb{E}[P_{12}|\mathcal{F}_{k+1, j}] & = \frac{1}{S_2^2} \sum\nolimits_{i=1}^{S_2}\mathbb{E}[\| D_{x}(u_{k, j}; \bar{\xi}_i^{x}) - D_{x}(u_{k, j}) \|^2 |\mathcal{F}_{k+1, j}]. \label{eq42-2}
\end{align}

Notice that from step \ref{line: 3-2} of Algorithm \ref{alg:3}, and by the non-expansive property of projection, it follows that $\|v_{k, j} - v_{k, j-1}\| \le \lambda_2 \eta \|D_{k, j-1}^{v}\|$. Using this inequality, (\ref{eq11}), Assumption \ref{assum:4}, $\|v_{k, j-1}\|\le M$, and $\bar{M}_1$ in Lemma \ref{lemma:3}, we get, for any $i$,
\begin{align}\label{eq43}
&\|D_{x}(u_{k, j}; \bar{\xi}_i^{x}) - D_{x}(u_{k, j-1};\bar{\xi}_i^{x})\|^2 \\ \nonumber
&\le 2\| \nabla_{x}F(x_{k, j}, y_{k, j}; \xi_i^{x}) -  \nabla_{x}F(x_{k, j-1}, y_{k, j-1}; \xi_i^{x})\|^2 \\ \nonumber
&~~~+ 4 \|\nabla_{x}\nabla_{y}G(x_{k, j}, y_{k, j}; \delta_{i}^{x})\|^2 \|v_{k, j} - v_{k, j-1}\|^2\\ \nonumber
& ~~~+ 4M^2 \|\nabla_{x}\nabla_{y}G(x_{k, j}, y_{k, j}; \delta_i^{x}) -  \nabla_{x}\nabla_{y}G(x_{k, j-1}, y_{k, j-1}; \delta_i^{x}) \|^2 \\ \nonumber
&\le \bar{M}_1(\|x_{k, j} - x_{k, j-1}\|^2 + \|y_{k, j} - y_{k, j-1}\|^2) +  \eta_{v}\|D_{k, j-1}^{v}\|^2.
\end{align}

Moreover, from (\ref{eq6}), (\ref{eq11}), and the fact that $\|v_{k, j}\|\le M$, it holds that for any $i$, 
\begin{align}\label{eq43-1}
&\mathbb{E}[\|D_{x}(u_{k, j}; \bar{\xi}_i^{x}) - D_{x}(u_{k, j})\|^2 | \mathcal{F}_{k+1, j}] \\ \nonumber
& \le 2 \mathbb{E}[\|\nabla_{x}F(x_{k,j}, y_{k, j}; \xi_i^{x})-  \nabla_{x}f(x_{k, j}, y_{k, j})\|^2 |\mathcal{F}_{k+1, j}] \\ \nonumber
& ~~~+ 2M^2 \mathbb{E}[\|\nabla_{x}\nabla_{y}G(x_{k, j}, y_{k, j}; \delta_i^{x})- \nabla_{x}\nabla_{y}g(x_{k, j}, y_{k, j}) \|^2|\mathcal{F}_{k+1, j}] \\ \nonumber
& \le 2 \sigma_f^2 + 2M^2\sigma_{gxy}^2.
\end{align}

Define $\sigma_1 = 4 (\tau^x)^2{S_2}^{-1}(\sigma_f^2 + M^2\sigma_{gxy}^2)$. Combining (\ref{eq41}), (\ref{eq42-0}), (\ref{eq42-1}), (\ref{eq42-2}), (\ref{eq43}), and (\ref{eq43-1}), and taking a full expectation, we obtain, for any non-negative integers $k$ and $j$, 
\begin{align}\label{eq44}
& \mathbb{E}\|D_{k, j}^{x} - D_{x}(u_{k, j})\|^2 \\  \nonumber
&\overset{(i)}{\le} (1 - \tau^x) \mathbb{E}\|D_{k, j-1}^{x} - D_{x}(u_{k, j-1}) \|^2 + 2{\eta_{v}}{S_2}^{-1}\mathbb{E}\|D_{k, j-1}^{v}\|^2 + \sigma_1 \\ \nonumber
& ~~~~+ 2{\bar{M}_1}{S_2}^{-1}(\mathbb{E}\|x_{k, j} - x_{k, j-1}\|^2 + \mathbb{E}\|y_{k, j} - y_{k, j-1}\|^2) \\ \nonumber
&\le (1 - \tau^x) \mathbb{E}\|D_{k}^{x} - D_{x}(u_{k})\|^2 + 2\eta_{v} {S_2}^{-1}\sum\nolimits_{i=0}^{j}\mathbb{E}\|D_{k, i-1}^{v}\|^2 + (j+1)\sigma_1 \\ \nonumber
&~~~~+ 2{\bar{M}_1}{S_2}^{-1}(\mathbb{E}\|x_{k+1} - x_k\|^2 + \mathbb{E}\|y_{k+1} - y_k\|^2)
\end{align}
where the last inequality follows from recursively applying $(i)$, under $D_{k, -1}^{x} = D_k^{x}$, $u_{k, -1} = u_k$, and $x_{k, j} = x_{k+1}$, $y_{k, j} = y_{k+1}$ for any non-negative integer $j$.

Moreover, due to $y_{k+1} = \bar{y}_{k,T-1}$, $y_{k} = \bar{y}_{k,-1}$, and step \ref{line:2-2} of Algorithm \ref{alg:2}, we observe that
\begin{align}\label{eq44-1}
\|y_{k+1} - y_{k}\|^2 \le T\sum\nolimits_{t=0}^{T-1}\|\bar{y}_{k, t}- \bar{y}_{k, t-1}\|^2  = T\lambda_1^2\beta^2 \sum\nolimits_{t=0}^{T-1}\|D_{k, t-1}^{y}\|^2.
\end{align}

Utilizing (\ref{eq44-1}) and step \ref{line:8} of Algorithm \ref{alg:1}, (\ref{eq44}) implies that, for any non-negative integers $k$ and $j$, 
\begin{align}\label{eq45}
\mathbb{E}\|D_{k, j}^{x} - D_{x}(u_{k, j})\|^2 & \le (1 - \tau^x) \mathbb{E}\|D_{k}^{x} - D_{x}(u_{k})\|^2 + \frac{2}{S_2}\eta_{v} \sum\nolimits_{i=0}^{j}\mathbb{E}\|D_{k, i-1}^{v}\|^2\\ \nonumber
&~~~~ + 2{S_2}^{-1}(\bar{M}_1\alpha^2 \mathbb{E}\|D_k^{x}\|^2 + \eta_{y}W_{k, T}^{y}) + (j+1)\sigma_1.
\end{align}

Then, $V_k^{x}$ can be bounded as follows. On one hand, if $\mod(k, q_1)=0$, by step \ref{line:4} of Algorithm \ref{alg:1}, (\ref{eq6}), the mutual independence of $\mathcal{S}_{1, x}$,
\begin{align}\label{eq46}
\mathbb{E}[\|D_k^{x} - D_{x}(u_k)\|^2| \mathcal{F}_{k}^{'}] &= {S_1}^{-2}\sum\nolimits_{i=1}^{S_1}\mathbb{E}[\|D_{x}(u_k; \bar{\xi}_i^{x}) - D_{x}(u_k)\|^2| \mathcal{F}_{k}^{'}].
\end{align}
Following a similar derivation as in (\ref{eq43-1}), and taking a full expectation, we get from (\ref{eq46}) that 
\begin{equation} \label{eq47}
V_k^{x} \le 2{S_1}^{-1}\sigma_f^2 + 2M^2{S_1}^{-1}\sigma_{gxy}^2.
\end{equation}

On the other hand, if $\mod(k, q_1)\neq 0$, we derive from line \ref{line:5} of Algorithm \ref{alg:1}, $u_{k-1, J-1} = u_k$, and (\ref{eq45}) that 
\begin{align}  \label{eq47-1}
V_k^{x} & = \mathbb{E}\|D_{k-1,J-1}^{x} - D_{x}(u_{k-1, J-1})\|^2  \le (1 - \tau^x) V_{k-1}^{x} + {2}/{S_2} W_{k-1, T,J} + J\sigma_1.
\end{align}
By using induction on this inequality, and leveraging the definition of $\sigma_1$ and (\ref{eq47}),  it can be concluded that
\begin{align} \label{eq48}
V_k^{x} & \le V_{(n_k-1)q_1}^{x} + C_0 + (\tau^x)^{-1}J\sigma_1\\ \nonumber
& \le \big({2}{S_1}^{-1} + 4{S_2}^{-1} \tau^{x} \big) J^2 (\sigma_f^2 + M^2\sigma_{gxy}^2) + C_0
\end{align}
where $C_0 =  \frac{2}{S_2}\sum\nolimits_{i=(n_k-1)q_1}^{k-1}W_{i, T,J}$, and $n_k$ is the integer such that $(n_k-1)q_1< k<(n_k)q_1$.

Since $\sum\nolimits_{i=s}^{s-1}(\cdot)=0$ for any non-negative integer $s$, (\ref{eq48}) also holds for the case $\mod(k, q_1)=0$ if $n_k$ is chosen to be the integer such that $k=(n_k-1)q_1$. Consequently, the inequality for $V_k^{x}$ is obtained. 

Second, we prove the inequality for $V_{k, T}^{y}$. From step \ref{line:2-4} of Algorithm \ref{alg:2}, for any non-negative integers $k$ and $t$, 
\begin{align}\label{eq49}
&\mathbb{E}[\|D_{k, t}^{y} - D_{y}(\bar{u}_{k, t})\|^2| \mathcal{F}_{k+1, t}] \\ \nonumber
& = (1 - \tau^y)^2 \|D_{k, t-1}^{y} - D_{y}(\bar{u}_{k, t-1})\|^2  +  \mathbb{E}[P_2|\mathcal{F}_{k+1, t}]
\end{align}
with $P_2 = \| D_{y}(\bar{u}_{k,t}; \mathcal{S}_{2, y})- D_{y}(\bar{u}_{k,t}) + (1 -\tau^y) \big(D_{y}(\bar{u}_{k,t-1}) - D_{y}(\bar{u}_{k,t-1}; \mathcal{S}_{2, y})\big)\|^2$. 

Define $\sigma_2 = 2(\tau^y)^2 S_2^{-1}\sigma_g^2$. Following derivations similar to (\ref{eq42-0}), (\ref{eq42-1}), (\ref{eq42-2}), (\ref{eq43}), (\ref{eq43-1}), and (\ref{eq44}), we obtain from (\ref{eq49}) that, for any non-negative integers $k$ and $t$,
\begin{align}\label{eq51}
\mathbb{E}\|D_{k, t}^{y} - D_{y}(\bar{u}_{k, t})\|^2 &\le (1 - \tau^y) \mathbb{E}\|{D}_{k, t-1}^{y} - D_{y}(\bar{u}_{k, t-1})\|^2 + \sigma_2 \\ \nonumber
& ~~~+ {2}{S_2}^{-1}L_g^2(\mathbb{E}\|\bar{x}_{k, t}- \bar{x}_{k, t-1}\|^2 + \lambda_1^2 \beta^2\mathbb{E}\|D_{k, t-1}^{y}\|^2)\\ \nonumber
&\le (1 - \tau^y) \mathbb{E}\|D_k^{y} - D_{y}(x_k, y_k)\|^2 +  2L_g^2\alpha^2{S_2}^{-1}\mathbb{E}\|D_k^{x}\|^2 \\ \nonumber
&~~~+ 2 L_g^2\lambda_1^2\beta^2{S_2}^{-1}\sum\nolimits_{i=0}^{t}\mathbb{E}\|D_{k, i-1}^{y}\|^2 + (t+1) \sigma_2 
\end{align}
where the last inequality is derived by recursively applying the first inequality, combined with step \ref{line:8} of Algorithm \ref{alg:1}, $D_{k, -1}^{y}= D_k^{y}$, $\bar{u}_{k, -1}=(x_k, y_k)$, and $\bar{x}_{k, t}= x_{k+1}$ for any non-negative integer $t$.


Then, for any non-negative integer $k$ and positive integer $T$,  
\begin{align}\label{eq52}
V_{k, T}^{y} &= V_k^y + \sum\nolimits_{t=1}^{T-1}\mathbb{E}\|D_{k, t-1}^{y} - D_{y}(\bar{u}_{k, t-1})\|^2\\ \nonumber
&\le T V_k^y  + 2{S_2}^{-1}L_g^2(T-1) W_{k, T} + (T-1)T\sigma_2
\end{align} 
where the first equality results from the definition of $V_{k, T}^{y}$, (\ref{eq26-2}), and the fact that $D_{k, -1}^{y}= D_k^{y}$, $\bar{u}_{k, -1}=(x_k, y_k)$, the last inequality derives from (\ref{eq51}), $\sum\nolimits_{t=1}^{0}(\cdot)=0$, and $\sum\nolimits_{t=1}^{T-1}\sum\nolimits_{i=0}^{t-1}\mathbb{E}\|D_{k, i-1}^{y}\|^2 \le (T-1)\sum\nolimits_{t=0}^{T-1}\mathbb{E}\|D_{k, t-1}^{y}\|^2$.

In the following, we analyze the first term on the r.h.s. of (\ref{eq52}) through two cases. On one hand, if $\mod(k, q_1)=0$, with similar derivations as in (\ref{eq42-2}) and (\ref{eq43-1}), we get 
\begin{align}\label{eq52-1}
V_k^y  & \le {S_1}^{-1}\sigma_{g}^2.
\end{align}

On the other hand, if $\mod(k,q_1)\neq 0$, utilizing step \ref{line:5} of Algorithm \ref{alg:1}, (\ref{eq51}), and $\bar{u}_{k-1, T-1}=(x_k, y_k)$, it holds that
\begin{align}\label{eq53}
V_k^y &= \mathbb{E}\|D_{k-1, T-1}^{y} -  D_{y}(\bar{u}_{k-1, T-1}) \|^2 \\ \nonumber
&\le (1 - \tau^y) V_{k-1}^y + {2}{S_2}^{-1}L_g^2 W_{k-1, T} + T \sigma_2 .
\end{align} 
By using induction on (\ref{eq53}), and leveraging (\ref{eq52-1}) and the definition of $\sigma_2$, we get
\begin{align} \label{eq54}
V_k^y &\le V_{(n_k-1)q_1}^{y} + {2}{S_2}^{-1}L_g^2 \sum\nolimits_{i=(n_k-1)q_1}^{k-1}W_{i, T} + (\tau^y)^{-1}T\sigma_2  \\ \nonumber
&\le {2} {S_2}^{-1}L_g^2 \sum\nolimits_{i=(n_k-1)q_1}^{k-1}W_{i, T} + \left({S_1}^{-1} + {2}{S_2}^{-1}T \tau^y \right)\sigma_{g}^2
\end{align}
where $n_k$ is the integer such that $(n_k-1)q_1 < k < (n_k)q_1$.

Since $\sum\nolimits_{i=s}^{s-1}(\cdot)=0$ for any non-negative integer $s$, (\ref{eq54}) also holds for the case $\mod(k, q_1)=0$ if $n_k$ is chosen to be the integer such that $k=(n_k-1)q_1$. By combining (\ref{eq52}) and (\ref{eq54}), we obtain the inequality for $V_{k, T}^{y}$.

Next, we prove the inequality for $V_{k, J}^{v}$. From step \ref{line: 3-4} of Algorithm \ref{alg:3}, for any non-negative integers $k$ and $j$, 
\begin{align}\label{eq55}
&\mathbb{E}[\|D_{k, j}^{v} - D_{v}(u_{k, j})\|^2| \mathcal{F}_{k+1, j}] \\ \nonumber
&= (1 - \tau^v)^2 \|D_{k, j-1}^{v} - D_{v}(u_{k, j-1})\|^2 + \mathbb{E}[P_3|\mathcal{F}_{k+1, j}]
\end{align}
with $P_3 = \|D_{v}(u_{k, j}; \mathcal{S}_{2, v})- D_{v}(u_{k, j}) + (1 - \tau^v) \big(D_{v}(u_{k, j-1}) - D_{v}(u_{k, j-1}; \mathcal{S}_{2, v})\big)\|^2$.

By similar derivations as in (\ref{eq42-0}), (\ref{eq42-1}), (\ref{eq42-2}), (\ref{eq43}), (\ref{eq43-1}), (\ref{eq44}), (\ref{eq44-1}), and (\ref{eq45}), the equality in (\ref{eq55}) implies that, for any non-negative integers $k$ and $j$, 
\begin{align} \label{eq57}
\mathbb{E}\|D_{k, j}^{v} - D_{v}(u_{k, j})\|^2 & \le (1 - \tau^v)\mathbb{E}\|D_{k}^{v} - D_{v}(u_{k})\|^2 + \frac{2}{S_2}\eta_{v}\sum\nolimits_{i=0}^{j}\mathbb{E}\|D_{k, i-1}^{v}\|^2\\ \nonumber
&~~~~ + 2{S_2}^{-1}(\bar{M}_1\alpha^2 \mathbb{E}\|D_k^{x}\|^2 + \eta_{y}W_{k, T}^{y}) + (j+1)\sigma_3
\end{align}
where $\sigma_3 = 4(\tau^v)^2 S_2^{-1} (\sigma_f^2 + M^2 \sigma_{gyy}^2)$.

Then, for any non-negative integer $k$, and positive integers $T$ and $J$,  
\begin{align}\label{eq58}
V_{k, J}^{v} & = V_k^v + \sum\nolimits_{j=1}^{J-1}\mathbb{E}\|D_{k, j-1}^{v}- D_{v}(u_{k, j-1})\|^2 \\ \nonumber
& \le J V_k^v + {2}{S_2}^{-1}(J-1) W_{k, T, J} + (J-1)J\sigma_3
\end{align}
where the first equality uses the definition of $V_{k,J}^{v}$, (\ref{eq26-2}), $D_{k, -1}^{v}= D_k^{v}$, $u_{k, -1}= u_k$, the last inequality results from  $\sum\nolimits_{j=1}^{J-1}\sum\nolimits_{i=0}^{j-1}\mathbb{E}\|D_{k,i-1}^{v}\|^2 \le (J-1)\sum\nolimits_{j=0}^{J-1}\mathbb{E}\|D_{k,j-1}^{v}\|^2$, $\sum\nolimits_{j=1}^{0}(\cdot)=0$, and (\ref{eq57}).

For the first term on the r.h.s. of (\ref{eq58}), with derivations analogous to (\ref{eq46}), (\ref{eq47}), (\ref{eq47-1}), and (\ref{eq48}), one has the following inequalities 
\begin{align}
&V_k^v  \le {2}{S_1}^{-1}(\sigma_f^2+ M^2 \sigma_{gyy}^2)\label{eq58-1} \\ 
&V_k^v \le (1 - \tau^v) V_{k-1}^v + {2}S_2^{-1} W_{k-1, T, J} + J \sigma_3, \label{eq58-2}\\ 
&V_k^v   \le {2}{S_2}^{-1}\sum\nolimits_{i=(n_k-1)q_1}^{k-1}W_{i,T, J}  + {2}{S_1}^{-1}(\sigma_f^2+ M^2 \sigma_{gyy}^2) + ({\tau^v})^{-1}J\sigma_3 \label{eq58-3}
\end{align}
where (\ref{eq58-1}) applies when $\mod(k,q_1) = 0$, while (\ref{eq58-2}) applies when $\mod(k,q_1) \neq 0$. (\ref{eq58-3}) holds for any non-negative integer $k$, with $n_k$ being the integer satisfying $(n_k-1)q_1 <k < n_k q_1$ if $\mod(k, q_1)\neq0$, and $k=(n_k-1)q_1$ if $\mod(k, q_1)=0$.

By (\ref{eq58}), (\ref{eq58-3}), and the definition of $\sigma_3$, the inequality for $V_{k, J}^{v}$ is obtained.
\end{proof}

\subsubsection{Convergence and Complexity Analysis}
The convergence analysis of ALS-SPIDER involves the recursive QP solution error $\delta_k^v$, which does not appear in the analysis of VRBO. To this end, for any non-negative integer $k$, we define the following Lyapunov function 
\begin{equation} \label{eq19}
H_k = \Phi(x_k) + \|y_k - y^*(x_k)\|^2 + \|v_k - v^*(x_k)\|^2
\end{equation}
as adopted by \cite{34}. Based on the above lemmas, the difference between $\mathbb{E}[H_{k+1}]$ and  $\mathbb{E}[H_{k}]$ over $k=0, \ldots, K-1$ can be obtained as follows.

\begin{lemma}\label{lemma:6}
Suppose Assumptions \ref{assum:1}, \ref{assum:2}, \ref{assum:3}, and \ref{assum:4} hold. Define $\Delta \tilde{H}_k = \mathbb{E}[H_{k+1}- H_k]$. Let $\beta$, $\eta$, $\lambda_1$, and  $\lambda_2$ satisfy the conditions in Lemma \ref{lemma:4}, $\tau^x$, $\tau^y$, and $\tau^v$ satisfy the conditions in Lemma \ref{lemma:5}, and $\sum_{i=s}^{s-1}(\cdot)=0$ for any non-negative integer $s$. Then, for any positive integers $K$, $T$, and $J$, 
\begin{align*}
&\sum\nolimits_{k=0}^{K-1}\left(\Delta \tilde{H}_k + \frac{\alpha}{2}\mathbb{E}\|\nabla \Phi(x_k)\|^2 - L_1 \mathbb{E}\|D_k^{x}\|^2 - L_2 W_{k, T}^{y} - L_3 W_{k, J}^{v} - L_4 \delta_k^{y} - L_5 \delta_k^{v}\right)\\ \nonumber
& \le \sum\nolimits_{k=0}^{K-1}\big[\alpha V_k^x + (1 + C_1(J))\beta_T^y V_{k, T}^y + \beta_J^v V_{k, J}^v\big]
\end{align*}
where $L_1 = -\frac{\alpha}{2}(1-\alpha L_{\Phi}) + \big(1 + C_1(J)\big)\beta_{yx} + \beta_{vx}$, 
\begin{align*}
&L_2= - (1/2)^{2T-1}\lambda_1^2\beta \left(1 + C_1(J)\right), \qquad \qquad ~~~\qquad L_3=  - (1/2)^{2J-1}\lambda_2^2\eta, \\ \nonumber
&L_4=\bar{M}_1\alpha + \big(1 + C_1(J)\big)\beta_{yy} - \big(1 - C_1(J)\big), \qquad ~~L_5 = 4L_g^2\alpha + \beta_{vv}-1.
\end{align*}
Here, $W_{k, T}^{y}$, $W_{k, J}^{v}$, $V_{k, T}^{y}$, and $V_{k, J}^{v}$ are in (\ref{eq26-2}) with $i=k$, $\delta_{k}^{y}$, $\delta_{k}^{v}$, $V_k^x$, and $\bar{M}_1$ are in Lemma \ref{lemma:3}, $\beta_{T}^{y}$, $\beta_{J}^{v}$, $\beta_{yx}$, $\beta_{yy}$, $\beta_{vx}$, $\beta_{vv}$, and $C_1(J)$ are in Lemma \ref{lemma:4}, and $L_{\Phi}$ is in Lemma \ref{lemma:2}.
\end{lemma}

\begin{proof}
Let positive integers $K$, $T$, and $J$ be arbitrary. Using the definition of $H_k$ in (\ref{eq19}), and the definitions of $\Delta \Phi_k$, $\delta_k^{y}$, and  $\delta_k^{v}$ in Lemma \ref{lemma:3}, it follows that
\begin{align}\label{eq60}
\sum\nolimits_{k=0}^{K-1}\mathbb{E}[H_{k+1} - H_{k}] = \sum\nolimits_{k=0}^{K-1}(\Delta \Phi_k + \delta_{k+1}^{y} -  \delta_{k}^{y}+ \delta_{k+1}^{v} -  \delta_{k}^{v}).
\end{align}

Combining (\ref{eq60}) and the inequality for $\delta_{k+1}^{v}$ in Lemma \ref{lemma:4}, we have, for any non-negative integer $k$, 
\begin{align}\label{eq61}
&\sum\nolimits_{k=0}^{K-1}\mathbb{E}[H_{k+1} - H_{k}] \le \sum\nolimits_{k=0}^{K-1}\bigg(\Delta \Phi_k + (1 + C_1(J))\delta_{k+1}^{y} - (1 - C_1(J))\delta_{k}^{y} \\ \nonumber
&\qquad \qquad \qquad \qquad \qquad \qquad \qquad ~~~~+ (\beta_{vv}-1)\delta_{k}^{v} + P_{4, k}\bigg)
\end{align}
where $P_{4, k} =\beta_{vx}\mathbb{E}\|D_k^{x}\|^2 - (\frac{1}{2})^{2J-1}\lambda_2^2 \eta W_{k, J}^{v}+ \beta_J^v V_{k, J}^v$.

Then, the proof is completed by combining (\ref{eq61}) and (\ref{eq22}), along with the inequality for $\delta_{k+1}^{y}$ from Lemma \ref{lemma:4}, and utilizing the definition of $P_{4, k}$.
\end{proof}

Next, we provide the convergence of ALS-SPIDER.
\begin{theorem}\label{tho:1}
Apply ALS-SPIDER to solve problem (\ref{eq1}). Suppose Assumptions \ref{assum:1}, \ref{assum:2}, \ref{assum:3}, and \ref{assum:4} hold. Let $S_2 = q_1$, $0< \beta <1$, $0< \eta <1$, $\tau^x = 0$, $\tau^y = 0$, $\tau^v = 0$, and $\sum_{i=s}^{s-1}(\cdot)=0$ for any non-negative integer $s$. Choose $\lambda_1 = \min\{\frac{1}{2L_g + 2\mu}, \frac{(1/2)^{2T+5}\mu r}{T L_g^2 \beta^2 C_3}\}$, $\lambda_2 = \min\left\{\frac{1}{2L_g + 2\mu}, \frac{\mu^2 \lambda_1 \beta}{128 \bar{M}_1 J \eta}, \frac{(1/2)^{2(T+J)+6}\mu r}{J(L_g^2 \eta^2 + T \bar{M}_1 \beta \eta)}\right\}$, and 
\begin{align}\label{eq62}
& \alpha = \min\bigg\{\frac{\mu \lambda_1 \beta}{96 C_{y}^2 + 8\bar{M}_1}, \frac{\mu \lambda_2 \eta}{48(L_g^2 + C_{v}^2)}, \frac{1}{6 L_{\Phi}}, \frac{\sqrt{r}}{12 \sqrt{\bar{M}_1}}, \frac{\mu}{C_4}, \frac{\mu}{C_5}, \frac{(1/2)^{2(T+J)+3}r}{\bar{M}_1 T \beta + L_g^2 \eta}\bigg\}, 
\end{align}
where $r = S_2/q_1$, $C_3 = J^2 + J^2 M^2 + 2T^2$, 
\begin{align} 
&C_4 = 192 \lambda_1 \beta L_g^2 ( 1 + {r}^{-1} C_3 ), ~~C_5 = 192 \lambda_2 \eta \bar{M}_1 (1 + C_y^2 + {r}^{-1}J ),  \label{eq64-2}
\end{align}
$M$ is in Definition \ref{definition:3}, $\bar{M}_1$ is in Lemma \ref{lemma:3}, and $L_{\Phi}$, $C_{y}$, and $C_{v}$ are given in Lemma \ref{lemma:2}. Then, for any given positive integers $T$ and $J$, 
\begin{align*}
\frac{1}{K}\sum\nolimits_{k=0}^{K-1}\mathbb{E}\|\nabla \Phi(x_k)\|^2 \le O\left(\frac{1}{K} + \frac{1}{S_1}\right).
\end{align*}
\end{theorem}

\begin{proof}
Replacing $V_k^x$, $V_{k, T}^y$, $V_{k, J}^v$ in Lemma \ref{lemma:6} with their upper bounds in (\ref{eq40-1}), (\ref{eq40-2}), and (\ref{eq40-3}), respectively, utilizing the inequalities $\sum\nolimits_{k=0}^{K-1}\sum\nolimits_{i=(n_k-1)q_1}^k W_{i, T} \le q_1 \sum\nolimits_{k=0}^{K-1}W_{k, T}$, $\sum\nolimits_{k=0}^{K-1}\sum\nolimits_{i=(n_k-1)q_1}^k W_{i, T, J} \le q_1 \sum\nolimits_{k=0}^{K-1}W_{k, T, J}$, and the definitions of $W_{k, T, J}$, $W_{k, T}$ and $\eta_y$, one has 
\begin{align}\label{eq64-3}
&\sum\nolimits_{k=0}^{K-1}\bigg(\Delta \tilde{H}_k + \frac{\alpha}{2}\mathbb{E}\|\nabla \Phi(x_k)\|^2 - \bar{L}_{1} \mathbb{E}\|D_k^{x}\|^2 - \bar{L}_{2} W_{k, T}^{y} - \bar{L}_{3} W_{k, J}^{v} - L_4 \delta_k^{y} - L_5 \delta_k^{v}\bigg)\\ \nonumber
& \le \left[\frac{2}{S_1}(\alpha + \beta_{T}^{y} + \beta_{J}^{v}) + \frac{8}{S_2}(\alpha \tau^x + \beta_{T}^{y}\tau^y + \beta_{J}^{v} \tau^v )\right] K \tilde{C}_1(J) (\sigma_f^2 + \sigma_g^2 + \sigma_{gxy}^2 + \sigma_{gyy}^2)
\end{align}
where $\bar{L}_{1} = L_1 + {2}{S_2}^{-1}q_1 \alpha^2 \left(\bar{M}_1 C_2(J) + \beta_T^y L_g^2 \tilde{C}_1(J)\right)$,
\begin{align}\label{eq64-4}
&\bar{L}_2 = L_2 + {2}{S_2}^{-1}q_1 T \lambda_1^2 \beta^2 \left(\bar{M}_1 C_2(J) + \beta_T^y L_g^2  \tilde{C}_1(J) \right), \qquad C_2(J) = \alpha + J \beta_J^v, \\ \nonumber
&\bar{L}_3 = L_3 + 2{S_2}^{-1}q_1\eta_v C_2(J), ~~~~\tilde{C}_1(J) = \max\left\{J^2(1 + M^2), T^2(1 + C_1(J))\right\}.
\end{align}

To establish the convergence of the algorithm, we analyze $\bar{L}_1$, $\bar{L}_2$, $\bar{L}_3$, $L_4$, and $L_5$. By the choices of $r_{y}$ and $r_{v}$ in Lemma \ref{lemma:4}, $\beta_{yy}$, $\beta_{yx}$, $\beta_{vv}$, and $\beta_{vx}$ in Lemma \ref{lemma:4} satisfy: $\beta_{yy} \le 1 - \frac{\mu}{4}\lambda_1 \beta$, $\beta_{yx} \le \frac{4}{\mu \lambda_1 \beta}C_{y}^2 \alpha^2 + L_g^2 \alpha^2 \beta_T^y$, $\beta_{vv} \le 1 - \frac{\mu}{4}\lambda_2 \eta$, $\beta_{vx} \le \frac{4}{\mu \lambda_2 \eta}C_v^2 \alpha^2 + \bar{M}_1 (1 + 2C_y^2)\alpha^2 \beta_J^v$. Using these inequalities, the definitions of $L_1$, $L_4$, and $L_5$ in Lemma \ref{lemma:6}, it follows that 
\begin{align*}
&L_1  \le  -\frac{\alpha}{2}(1 - \alpha L_{\Phi}) + \big(1 + C_1(J)\big) \frac{4}{\mu \lambda_1 \beta}C_y^2 \alpha^2 + \frac{4}{\mu \lambda_2 \eta}C_v^2 \alpha^2  \\ \nonumber
& \qquad ~~+ (1 + C_1(J))L_g^2 \alpha^2 \beta_T^y + \bar{M}_1(1 + 2C_y^2)\alpha^2 \beta_J^v\\ \nonumber
& L_4 \le \bar{M}_1 \alpha + 2C_1(J) - {\mu}\lambda_1 \beta \big(1 + C_1(J)\big)/4, \qquad  L_5 \le 4L_g^2 \alpha - {\mu}\lambda_2 \eta/4.
\end{align*}

Let
\begin{align}\label{eq65}
\alpha \le \min\bigg\{\frac{\mu \lambda_1 \beta}{96 C_{y}^2 + 8\bar{M}_1}, \frac{\mu \lambda_2 \eta}{48(L_g^2 + C_{v}^2)}, \frac{1}{6 L_{\Phi}} \bigg\},\qquad  \lambda_2 \le \frac{\mu^2 \lambda_1 \beta}{128\bar{M}_1J\eta}.
\end{align}
Together with the choices for $\lambda_1$ and $\beta$, the definition of $C_3$, $C_1(J)$ in Lemma \ref{lemma:4}, and $\tilde{C}_1(J)$ in (\ref{eq64-4}), it is easy to verify that 
\begin{align}
&C_1(J) \le \mu \lambda_1 \beta/16 < 1/16, ~~~ \tilde{C}_1(J) \le C_3,\label{eq65-1} \\ 
&L_1 \le -{\alpha}/{4} + 2 L_g^2 \alpha^2 \beta_T^y + \bar{M}_1(1 + 2C_y^2)\alpha^2 \beta_J^v, ~~~ L_4 \le 0, ~~~ L_5 \le 0.\label{eq65-2}
\end{align}

Moreover, using (\ref{eq65-1}), along with the definitions of $C_2(J)$ in (\ref{eq64-4}) and $\eta_v$ in Lemma \ref{lemma:5}, we derive 
\begin{align}
&\bar{M}_1 C_2(J) + \beta_T^y L_g^2 \tilde{C}_1(J) \le \bar{M}_1 (\alpha + J\beta_J^v) + \beta_T^y L_g^2 C_3, \label{eq65-3}\\ 
&\eta_v C_2(J) = 4L_g^2 \lambda_2^2 \eta^2 (\alpha + J\beta_J^v). \label{eq65-4}
\end{align}

By (\ref{eq65-1}), (\ref{eq65-2}), (\ref{eq65-3}), (\ref{eq65-4}), the definitions of $L_2$ and $L_3$ in Lemma \ref{lemma:6}, $\beta_T^y$ and $\beta_J^v$ in Lemma \ref{lemma:4}, and $r$ in Theorem \ref{tho:1}, we have 
\begin{align*}
&\bar{L}_1 \le -\frac{\alpha}{4} + 2r^{-1}\bar{M}_1 \alpha^3 + {16}{\mu}^{-1}\big[L_g^2 \lambda_1 \beta(1 + r^{-1}C_3) + \bar{M}_1 \lambda_2 \eta(1 + C_y^2 + r^{-1}J)\big]\alpha^2, \\ \nonumber
& \bar{L}_2 \le \frac{2}{r}\left(\bar{M}_1 \alpha + \frac{8}{\mu}\lambda_2 \eta \bar{M}_1 J  + \frac{8}{\mu}\lambda_1 \beta L_g^2 C_3 \right) T \lambda_1^2 \beta^2   - ({1}/{2})^{2T-1}\lambda_1^2 \beta, \\ \nonumber
& \bar{L}_3 = \frac{8}{r} \bigg(\alpha + \frac{8}{\mu}\lambda_2 \eta J\bigg)L_g^2 \lambda_2^2 \eta^2 - ({1}/{2})^{2J-1}\lambda_2^2 \eta.
\end{align*}

Let 
\begin{align}\label{eq66}
\alpha \le \min\bigg\{\frac{\mu}{C_4}, \frac{\mu}{C_5}, \bar{\alpha}_1, \bar{\alpha}_2 \bigg\},  \lambda_1 \le \frac{(1/2)^{2T+5}\mu r}{TL_g^2 \beta^2 C_3},  \lambda_2 \le \frac{(1/2)^{2(T+J)+6}\mu r}{J(L_g^2 \eta^2+ T\bar{M}_1 \beta \eta)}
\end{align}
where $\bar{\alpha}_1 = \frac{\sqrt{r}}{12 \sqrt{\bar{M}_1}}$, $\bar{\alpha}_2 = \frac{(1/2)^{2(T+J)+3}r}{\bar{M}_1 T \beta + L_g^2 \eta}$, and $C_4$ and $C_5$ are in (\ref{eq64-2}). One can see that $\bar{L}_1 \le 0$, $\bar{L}_2 \le 0$, $\bar{L}_3 \le 0$.

It is obvious that under the settings of Theorem \ref{tho:1}, $\alpha$, $\lambda_1$, and $\lambda_2$ satisfy the conditions specified in (\ref{eq65}) and (\ref{eq66}). With $\bar{L}_i\le 0$ for $i=1$, $2$, $3$, $L_j\le 0$ for $j=4$, $5$, and $\tau^x$, $\tau^y$, and $\tau^v$ set to $0$, we get, from (\ref{eq64-3}), 
\begin{align*}
\sum\nolimits_{k=0}^{K-1}\big(\Delta \tilde{H}_k + \frac{\alpha}{2}\mathbb{E}\|\nabla \Phi(x_k)\|^2\big) \le \frac{2K}{S_1}\tilde{C}_1(J)(\alpha + \beta_T^y + \beta_J^v)(\sigma_f^2 + \sigma_g^2 +  \sigma_{gxy}^2 + \sigma_{gyy}^2).
\end{align*}
Notice that by the definition of $\Delta \tilde{H}_k$ and the inequality $H_K \ge \inf\nolimits_{x \in \mathbb{R}^p}\Phi(x)$, it holds that 
$-\sum\nolimits_{k=0}^{K-1}\Delta \tilde{H}_k \le H_0 - \inf\nolimits_{x\in \mathbb{R}^p}\Phi(x)$. Substituting this relation into the preceding inequality, and using (\ref{eq65-1}) and the definitions of $\beta_T^y$ and $\beta_J^v$, it holds that
\begin{align}\label{eq67}
\frac{\alpha}{2}\sum\nolimits_{k=0}^{K-1}\mathbb{E}\|\nabla \Phi(x_k)\|^2 \le  H_0 - \inf\limits_{x \in \mathbb{R}^p}\Phi(x) + {2K}{S_1}^{-1}P_5
\end{align}
where $P_5 = C_3(\alpha + \frac{8}{\mu}\lambda_1 \beta + \frac{8}{\mu}\lambda_2 \eta)(\sigma_f^2 + \sigma_g^2 +  \sigma_{gxy}^2 + \sigma_{gyy}^2)$. By multiplying both sides of (\ref{eq67}) by $\frac{2}{K\alpha}$, we obtain
\begin{align*}
\frac{1}{K}\sum\nolimits_{k=0}^{K-1}\mathbb{E}\|\nabla \Phi(x_k)\|^2 \le \frac{2}{K\alpha}\big(H_0 - \inf\limits_{x\in \mathbb{R}^p}\Phi(x)\big) + \frac{4}{\alpha S_1}P_5
\end{align*}
which implies that, for any given positive integers $T$ and $J$, $\frac{1}{K}\sum\nolimits_{k=0}^{K-1}\mathbb{E}\|\nabla \Phi(x_k)\|^2 \le O(\frac{1}{K}+ \frac{1}{S_1})$.
Then, the proof is completed.
\end{proof}

Theorem \ref{tho:1} shows that for any given positive integers $T$, $J$, ALS-SPIDER converges sublinearly w.r.t. the number of outer loop iterations $K$, and the convergence error decays sublinearly w.r.t. the outer batch size $S_1$. Below, we present the sample complexity for ALS-SPIDER to converge to an $\epsilon$-stationary point of problem (\ref{eq1}), by properly choosing hyperparameters.

\begin{corollary}\label{cor:1}
Apply ALS-SPIDER to solve problem (\ref{eq1}) under the same assumptions and parameter settings as in Theorem \ref{tho:1}. Choose $K = O(\epsilon^{-1})$, $S_1 = O(\epsilon^{-1})$, and $S_2 = q_1 = S_1^{1/2}(= O(\epsilon^{-0.5}))$. Then, for any positive integers $T$ and $J$, ALS-SPIDER finds an $\epsilon$-stationary point of problem (\ref{eq1}) with a sample complexity of $O(\epsilon^{-1.5})$.
\end{corollary}

\begin{proof}
From Theorem \ref{tho:1}, we have $\frac{1}{K}\sum\nolimits_{k=0}^{K-1}\mathbb{E}\|\nabla \Phi(x_k)\|^2 \le O\left(\frac{1}{K} + \frac{1}{S_1}\right)$. By setting $K = O(\epsilon^{-1})$ and $S_1 = O(\epsilon^{-1})$, for any given positive integers $T$ and $J$, ALS-SPIDER can find an $\epsilon$-stationary point for problem (\ref{eq1}), i.e., $\mathbb{E}\|\nabla \Phi(\bar{x})\|^2 \le \epsilon$, with $\bar{x}$ chosen uniformly from $\{x_0, \ldots, x_{K-1}\}$. Furthermore, the complexity of computing gradients $\nabla_{x}F(x, y; \xi)$ and $\nabla_{y}F(x, y; \xi)$ is given by Gc($F$, $\epsilon$)$ = 2 KS_1/q_1  + 4KS_2J = O(\epsilon^{-1.5})$, the complexity of computing the gradient $\nabla_{y}G(x, y; \delta)$ is given by Gc($G$, $\epsilon$)$ =  K S_1/q_1 + 2KS_2 T= O(\epsilon^{-1.5})$, the complexity of computing the Jacobian-vector product $\nabla_{x}\nabla_{y}G(x, y; \delta)v$ is given by  Jv($G$, $\epsilon$) $= KS_1/q_1 + 2KS_2J = O(\epsilon^{-1.5})$, and the complexity of computing the Hessian-vector product $\nabla_{y}^2G(x, y; \delta)v$ is given by  Hv($G$, $\epsilon$) $ = KS_1/q_1  + 2KS_2J = O(\epsilon^{-1.5})$. Therefore, the total sample complexity for ALS-SPIDER to reach an $\epsilon$-stationary point of problem (\ref{eq1}) is $O(\epsilon^{-1.5})$.
\end{proof}

The complexity result in Corollary \ref{cor:1} improves upon that of VRBO under the same assumptions and the same order for the batch sizes $S_1$, $S_2$, and the period $q_1$. This improvement is achieved by eliminating the requirement that the iteration step $J$ for estimating $\nabla \Phi(x)$ be $O(\log{\epsilon^{-1}})$. Moreover, our result matches the sample complexity of SPIDER \cite{51} for single-level nonconvex optimization, demonstrating that multiple LL updates do not theoretically reduce the efficiency of the algorithm for nonconvex-strongly-convex (n.s.c.) bilevel optimiztion. Additionally, it is worth noting that the conclusion in Corollary \ref{cor:1} holds independently of the iteration steps $T$ and $J$, providing flexibility in its implementation. To the best of our knowledge, this is the first stochastic algorithm to achieve such complexity for n.s.c. bilevel optimization when updating the LL variable multiple times. 

\section{A STORM-based alternating stochastic variance-reduced algorithm} \label{sec:4}

ALS-SPIDER achieves the optimal complexity of $O(\epsilon^{-1.5})$ under Assumptions \ref{assum:1}, \ref{assum:2}, \ref{assum:3}, and \ref{assum:4}. However, its convergence analysis requires a batch size $S_1$ of the order $O(\epsilon^{-1})$ for the periodic gradient estimate, and a batch size $S_2$ of the order $O(\epsilon^{-0.5})$ for Algorithms \ref{alg:2} and \ref{alg:3}. To avoid using large batches in each iteration, we employ the STORM variance reduction technique for gradient estimate instead of SPIDER. This is motived by the fact that in single-level optimization, STORM has the same complexity as SPIDER but imposes a weaker limitation on batch size \cite{50}. As a result, we obtain ALS-STORM by modifying the settings of $q_1$, $\tau^x$, $\tau^y$, and $\tau^v$ in Algorithm \ref{alg:1}, as shown in line \ref{line:1} of Algorithm \ref{alg:4}.

\begin{remark} \label{remark:3} 
(i) Since $q_1$ is set to $K$, ALS-STORM computes the gradient estimate in line \ref{line:4} of Algorithm \ref{alg:1} only when $k=0$. (ii) As mentioned in Section \ref{sec:1}, BSVRB$^{\text{v}2}$ \cite{61} employs the STORM variance reduction technique to perform one LL update and one hypergradient estimation iteration, and is shown to have a complexity of $O(\epsilon^{-1.5})$ in solving problem (\ref{eq1}). However, extending this algorithm to perform multiple LL updates and multiple hypergradient estimation iterations is not  straightforward. Indeed, unlike BSVRB$^{\text{v}2}$ which directly uses the gradient estiamtes from two consecutive iterations to construct variance-reduced estimates, ALS-STORM indirectly realizes variance reduction between two consecutive iterations by recursively applying variance reduction in multiple updates of LL or auxiliary variables. For example, the variance reduction between $V_k^y$ and $V_{k-1}^y$ in (\ref{eq53}) at the outer loop is recursively achieved through the inner loop variance reduction in (\ref{eq51}).
\end{remark}
\subsection{Convergence and Complexity Analysis}

A straightforward approach to establish the convergence of ALS-STORM is to extend the theoretical framework in Theorem \ref{tho:1} to the case where $\tau^x >0$, $\tau^y >0$, $\tau^v >0$. However, since $q_1$ for ALS-STORM is set to $K$, a constant batch size for $S_2$ will result in the choice of $\alpha$ in (\ref{eq62}) being of the order $O(\frac{1}{K})$. As a result, this leads to a sample complxity of $O(\epsilon^{-2})$, which is worse than that of ALS-SPIDER. To address this, for any non-negative integer $k$, with $u_k = (x_k, y_k,v_k)$, we introdue the following Lyapunov function
\begin{equation} \label{eq68}
H_{1, k} = H_k + \frac{C_{\eta}^x}{\alpha} \|D_k^x - D_x(u_k)\|^2 + \frac{C_{\eta}^y}{\alpha} \|D_k^y - D_y(x_k, y_k)\|^2 +  \frac{C_{\eta}^v}{\alpha} \|D_k^v - D_v(u_k)\|^2.
\end{equation}

\begin{theorem}\label{tho:2}
Apply ALS-STORM to solve problem (\ref{eq1}). Suppose Assumptions \ref{assum:1}, \ref{assum:2}, \ref{assum:3}, and \ref{assum:4} hold. Let $q_1 = K$, $0 < \lambda_1 \le 1/(2L_g + 2\mu)$, $0 < \lambda_2 \le 1/(2L_g + 2\mu)$, $\beta = C_{\beta} \alpha$, $\eta = C_{\eta} \alpha$, $\tau^x = C_{\eta x}\alpha^2$, $\tau^y = C_{\eta y}\alpha^2$, $\tau^v = C_{\eta v}\alpha^2$,
\begin{align*}
C_{\eta}^x = \frac{1}{C_{\eta x}},~~~ C_{\eta}^y = \frac{16}{C_{\eta y}\mu}\lambda_1 C_{\beta}T, ~~~C_{\eta}^v = \frac{8}{C_{\eta v}\mu} \lambda_2 C_{\eta}J
\end{align*}
and $\sum_{i=s}^{s-1}(\cdot)=0$ for any non-negative integer $s$. Choose $S_1 = O(\epsilon^{-0.5})$, $C_{\eta} = \frac{48(L_g^2 + C_v^2)}{\mu \lambda_2}$, $C_{\beta} = \max \{ \frac{96C_y^2 + 8\bar{M}_1}{\mu \lambda_1}, \frac{128 \bar{M}_1 J C_{\eta} \lambda_2}{\mu^2 \lambda_1}\}$, $C_{\eta x} = \frac{1} {S_2} C_7$, $C_{\eta y} = \frac{1} {\mu S_2}C_7 \lambda_1 C_{\beta}T$, $C_{\eta v} = \frac{1} {\mu S_2} C_7 \lambda_2 C_{\eta}J$, and 
\begin{align*}
&\alpha = \min \bigg\{\sqrt{\frac{S_2 \epsilon}{16 \tilde{P}_6}}, \sqrt{\frac{(1/2)^{2(T+J+4)} \mu S_2}{\lambda_2 L_g^2 C_{\eta}^2 J + TC_{\beta}C_3 C_6}}, \frac{1}{32}\sqrt{\frac{\mu}{C_6(S_2^{-1}C_3 + 1 + C_y^2)}}, \sqrt{\frac{S_2}{C_7}}, \frac{1}{C_8}\bigg\}, 
\end{align*}
where $C_6 = L_g^2 \lambda_1 C_{\beta} + \bar{M}_1 \lambda_2 C_{\eta}$, $C_7 = 2^{2(T+J)+11}\left[L_g^2 C_{\eta} + (\bar{M}_1 + L_g^2)(1 + TC_{\beta})\right]$,  
\begin{align*} 
C_8 = 1 + 6 L_{\Phi} + C_{\beta}+ C_{\eta}+ \sqrt{C_{\eta x}} + \sqrt{C_{\eta y}}+ \sqrt{C_{\eta v}}, 
\end{align*}
$\tilde{P}_6$ is in (\ref{eq72-3}), $\bar{M}_1$ is in Lemma \ref{lemma:3}, $L_{\Phi}$, $C_{y}$, and $C_{v}$ are given in Lemma \ref{lemma:2}, and $C_3$ is in Theorem \ref{tho:1}. Then, for any given batch size $S_2$ and positive integers $T$ and $J$, ALS-STORM achieves an $\epsilon$-stationary point for problem (\ref{eq1}) with a sample complexity of $O(\epsilon^{-1.5})$.
\end{theorem}

\begin{proof}
With $\alpha \le {1}/{C_8}$, $0< \lambda_1 \le 1/(2L_g + 2\mu)$, and $0< \lambda_2 \le 1/(2L_g + 2\mu)$, the proofs in Lemmas \ref{lemma:3}, \ref{lemma:4}, \ref{lemma:5}, and \ref{lemma:6} remain valid for ALS-STROM.

Then, based on the definitions of $H_{1, k}$ in (\ref{eq68}), $\Delta \tilde{H}_k$ in Lemma \ref{lemma:6}, $V_k^x$ in Lemma \ref{lemma:3}, and $V_k^y$ and $V_k^v$ in (\ref{eq26-2}), we have 
\begin{align} \label{eq70}
\sum\nolimits_{k=0}^{K-1}\mathbb{E}[H_{1, k+1} - H_{1, k}] = \sum\nolimits_{k=0}^{K-1}\left(\Delta \tilde{H}_k + \Delta \tilde{V}_{k}^x + \Delta \tilde{V}_{k}^y + \Delta \tilde{V}_{k}^v\right)
\end{align}
where $\Delta \tilde{V}_{k}^x = \frac{C_{\eta}^x}{\alpha}(V_{k+1}^x - V_k^x)$, $\Delta \tilde{V}_{k}^y = \frac{C_{\eta}^y}{\alpha}(V_{k+1}^y - V_k^y)$, $\Delta \tilde{V}_{k}^v = \frac{C_{\eta}^v}{\alpha}(V_{k+1}^v - V_k^v)$.

Combining (\ref{eq70}) and the inequality for $\Delta \tilde{H}_k$ in Lemma \ref{lemma:6}, we obtain
\begin{align} \label{eq71}
&\sum\nolimits_{k=0}^{K-1}\mathbb{E}[H_{1, k+1} - H_{1, k}] \le \sum\nolimits_{k=0}^{K-1}\left[P_{5, k} + \alpha V_k^x + (1 + C_1(J))\beta_T^y V_{k, T}^y  \right. \\ \nonumber
&\left. \qquad \qquad \qquad \qquad \qquad \qquad \qquad \qquad + \beta_J^v V_{k, J}^v + \Delta \tilde{V}_{k}^x + \Delta \tilde{V}_{k}^y + \Delta \tilde{V}_{k}^v\right]
\end{align}
where $P_{5, k} = - \frac{\alpha}{2}\mathbb{E}\|\nabla \Phi(x_k)\|^2 + L_1 \mathbb{E}\|D_k^{x}\|^2 + L_2 W_{k, T}^{y} + L_3 W_{k, J}^{v} + L_4 \delta_k^{y} + L_5 \delta_k^{v}$.

Define $\Delta \tilde{H}_{1, k} = \mathbb{E}[H_{1, k+1} - H_{1,k}]$. By (\ref{eq71}), (\ref{eq47-1}), (\ref{eq52}),(\ref{eq53}), (\ref{eq58}), and (\ref{eq58-2}), and the definitions of $\sigma_1$, $\sigma_2$, $\sigma_3$, $W_{k, T,J}$, $W_{k,T}$, and $\eta_y$, it holds that 
\begin{align}\label{eq71-1}
&\sum\nolimits_{k=0}^{K-1}\big(\Delta \tilde{H}_{1, k} + \frac{\alpha}{2}\mathbb{E}\|\nabla \Phi(x_k)\|^2 - L_1^* \mathbb{E}\|D_k^x\|^2 - L_2^* W_{k, T}^y - L_3^* W_{k, J}^v - L_4 \delta_k^y - L_5 \delta_k^v \\ \nonumber
& \qquad \qquad  -L_6^* V_k^x - L_7^* V_k^y  - L_8^* V_k^v  \big) \\ \nonumber
& \le \frac{4}{S_2}\big[(\tau^x)^2 + (\tau^y)^2 + (\tau^v)^2\big]\left(C_9 + \beta_T^y \tilde{C}_1(J) + C_2(J)\right)KC_3(\sigma_f^2 + \sigma_g^2 + \sigma_{gxy}^2 + \sigma_{gyy}^2)
\end{align}
where $L_1^* = L_1 + {2}{S_2}^{-1}\alpha^2 \left(\bar{M}_1 C_2(J) + \beta_T^y L_g^2 \tilde{C}_1(J) + C_9(\bar{M}_1 + L_g^2)\right)$, 
\begin{align*}
&L_2^* = L_2 + {2}{S_2}^{-1} T \lambda_1^2 \beta^2 \left( \bar{M}_1 C_2(J) + \beta_T^y L_g^2  \tilde{C}_1(J) + C_9(\bar{M}_1 + L_g^2 )\right), \\ \nonumber
&L_3^* = L_3 + 2{S_2}^{-1} \eta_v (C_2(J) + C_9), \qquad \qquad L_6^* = \alpha - C_{\eta}^x \tau^x/\alpha,\\ \nonumber
&L_7^* = \beta_T^y T(1 + C_1(J)) - C_{\eta}^y \tau^y/\alpha, \qquad \qquad L_8^* = \beta_J^v J - C_{\eta}^v \tau^v/\alpha.
\end{align*}
$C_9 = \frac{1}{\alpha}(C_{\eta}^x+ C_{\eta}^y + C_{\eta}^v)$, $\tilde{C}_1(J)$ and $C_2(J)$ are in (\ref{eq64-4}), and $C_3$ is in Theorem \ref{tho:1}.

Given $\alpha \le {1}/{C_8}$, and the choices of $C_{\beta}$ and $C_{\eta}$ in Theorem \ref{tho:2}, it follows that $\alpha$ and $\lambda_2$ satisfy the inequality in (\ref{eq65}). Using the conditions $\beta<1$ and $\lambda_1 \le 1/(2L_g + 2\mu)$, we further derive the inequalities in (\ref{eq65-1}), (\ref{eq65-2}), (\ref{eq65-3}), and (\ref{eq65-4}).

Moreover, due to the definitions of $C_{\eta}^x$, $C_{\eta}^y$, $C_{\eta}^v$, $\tau^x$, $\tau^y$, $\tau^v$, $C_7$, and $C_9$, and the choices of $C_{\eta x}$, $C_{\eta y}$, and $C_{\eta v}$, we get
\begin{align}\label{eq72}
C_9 = \frac{25}{\alpha C_7} S_2 \le \frac{1}{\alpha}P_6, ~C_{\eta}^x \tau^x = \alpha^2, ~C_{\eta}^y \tau^y = \frac{16}{\mu}\lambda_1 C_{\beta}T\alpha^2,~ C_{\eta}^v \tau^v = \frac{8}{\mu}\lambda_2 C_{\eta}J\alpha^2,
\end{align}
where 
\begin{equation}\label{eq72-1}
P_6 = \min \bigg\{\frac{S_2}{32(\bar{M}_1 + L_g^2)}, \frac{(1/2)^{2T+2}S_2}{T C_{\beta}(\bar{M}_1 + L_g^2)}, \frac{(1/2)^{2J-1} S_2}{24 L_g^2 C_{\eta}}\bigg\}. 
\end{equation}

By (\ref{eq65-1}), (\ref{eq65-2}), (\ref{eq65-3}), (\ref{eq65-4}), (\ref{eq72}), $L_2$ and $L_3$ in Lemma \ref{lemma:6}, $\beta_T^y$ and $\beta_J^v$ in Lemma \ref{lemma:4}, $\eta_v$ in Lemma \ref{lemma:5}, and the definitions of $\beta$ and $\eta$, one has 
\begin{align*}
&L_1^* \le -\frac{3 \alpha}{16} + \frac{2}{S_2}\bar{M}_1 \alpha^3 + {16}{\mu}^{-1}\bigg[L_g^2 \lambda_1 C_{\beta}\big(1 + \frac{1}{S_2}C_3\big) + \bar{M}_1 \lambda_2 C_{\eta}\big(1 + C_y^2 + \frac{1}{S_2}J \big)\bigg]\alpha^3, \\ \nonumber
&L_2^* \le \frac{2}{S_2}\left(\bar{M}_1 + \frac{8}{\mu}\lambda_2 C_{\eta} \bar{M}_1 J  + \frac{8}{\mu}\lambda_1 C_{\beta} L_g^2 C_3 \right) T \lambda_1^2 C_{\beta}^2 \alpha^3  - \frac{3}{4}({1}/{2})^{2T-1}\lambda_1^2 C_{\beta}\alpha, \\ \nonumber
&L_3^* \le \frac{8}{S_2} \bigg(1 + \frac{8}{\mu}\lambda_2 C_{\eta}J\bigg)L_g^2 \lambda_2^2 C_{\eta}^2 \alpha^3 - \frac{2}{3}({1}/{2})^{2J-1}\lambda_2^2 C_{\eta}\alpha, \\ \nonumber
& L_6^* \le 0, ~L_7^* \le 0,~ L_8^* \le 0. 
\end{align*}

Since 
\begin{align*}
\alpha \le \min \bigg\{\sqrt{\frac{S_2}{C_7}}, \frac{1}{32}\sqrt{\frac{\mu}{C_6(S_2^{-1}C_3 + 1 + C_y^2)}}, \sqrt{\frac{(1/2)^{2(T+J+4)} \mu S_2}{\lambda_2 L_g^2 C_{\eta}^2 J + TC_{\beta}C_3 C_6}}\bigg\}, 
\end{align*}
we conclude that $L_1^* \le 0$, $L_2^* \le 0$, and $L_3^* \le 0$. 

Define $\sigma^2 = \sigma_f^2 + \sigma_g^2 + \sigma_{gxy}^2 + \sigma_{gyy}^2$. Using $L_i^* \le 0$ for $i=1$, $2$, $3$, $6$, $7$, $8$, $L_j \le 0$ for $j=4$, $5$, along with  (\ref{eq65-1}), (\ref{eq72}), $\beta_T^y$ and $\beta_J^v$ in Lemma \ref{lemma:4}, $C_2(J)$ in (\ref{eq64-4}), the definitions of $\tau^x$, $\tau^y$, $\tau^v$, $\beta$, and $\eta$, and the choice of $\alpha$, (\ref{eq71-1}) implies that 
\begin{align} \label{eq72-2}
&\sum\nolimits_{k=0}^{K-1}\big(\Delta \tilde{H}_{1, k} + \frac{\alpha}{2}\mathbb{E}\|\nabla \Phi(x_k)\|^2 \big) \le \frac{4}{S_2}K\tilde{P}_6 \alpha^3
\end{align}
where 
\begin{align}\label{eq72-3}
\tilde{P}_6 = (C_{\eta x}^2 + C_{\eta y}^2 + C_{\eta v}^2)\bigg(1 + P_6 + \frac{8}{\mu}\lambda_1 C_{\beta}C_3 + \frac{8}{\mu}\lambda_2 C_{\eta}J\bigg)C_3 \sigma^2.
\end{align}

By the definition of $\Delta \tilde{H}_{1, k}$, (\ref{eq47}), (\ref{eq52-1}), and (\ref{eq58-1}), one can see that 
\begin{equation}\label{eq72-4}
- \sum\nolimits_{k=0}^{K-1} \Delta \tilde{H}_{1, k}\le H_0 -  \inf\limits_{x\in \mathbb{R}^p}\Phi(x) + \frac{2P_7}{\alpha S_1 }
\end{equation}
where $P_7 = (C_{\eta}^x + C_{\eta}^y + C_{\eta}^v)(1 + M^2)\sigma^2$. 

Combining (\ref{eq72-4}) and (\ref{eq72-2}), and multiplying both sides by $\frac{2}{K\alpha}$, we obtain
\begin{align}\label{eq72-5}
\frac{1}{K}\sum\limits_{k=0}^{K-1}\mathbb{E}\|\nabla \Phi(x_k)\|^2 \le \frac{2}{K\alpha}(H_0 - \inf\limits_{x\in \mathbb{R}^p}\Phi(x)) + \frac{4}{K\alpha^2 S_1}P_7 + \frac{8}{S_2}\tilde{P}_6 \alpha^2.
\end{align}

Based on (\ref{eq72-5}) and the choices of $\alpha$ and $S_1$, to find an $\epsilon$-stationary point for problem (\ref{eq1}), it suffices to choose $K= O(\epsilon^{-1.5})$. Consequently, following a similar derivation as in Corollary \ref{cor:1}, the sample complexity of ALS-STORM is $O(\epsilon^{-1.5})$. Then, the proof is completed.
\end{proof}

Theorem \ref{tho:2} demonstrates that ALS-STORM has the same sample complexity as ALS-SPIDER. Notably, the convergence and complexity analysis of ALS-STORM only requires using  a large batch size at the intial iteration, making it more practical to implement compared to ALS-SPIDER.

\begin{algorithm}[h]   
\caption{ALS-STORM algorithm for problem (\ref{eq1})} \label{alg:4}
\begin{algorithmic}[1]
\STATE Set $q_1 = K$, $0 < \tau^{x} < 1$, $0 < \tau^{y} < 1$, and $0 < \tau^{v} < 1$ in Algorithm \ref{alg:1} \label{line:1}
\end{algorithmic}
\end{algorithm}

\section{Experiments}\label{sec:5}
In this section, we evaluate the performance of our proposed algorithms by comparing ALS-SPIDER with VRBO \cite{21}, and ALS-STORM with ALS-SPIDER. 

For clarity, we denote $q_1$ as the period, $S_1$ as the outer batch size, $S_2$ as the inner batch size, and $T$ and $J$ as the iteration steps for solving the LL problem and estimating $\nabla \Phi(x)$, respectively. Moreover, we use $\alpha$, $\beta$ and $\eta$ to represent the step sizes, while $\tau^{x}$, $\tau^{y}$ and $\tau^{v}$ denote the momentum parameters for variance-reduced estimates. Specifically, ($\alpha$, $\tau^{x}$), ($\beta$, $\tau^{y}$), and ($\eta$, $\tau^{v}$) are respectively used for the UL variable, LL variable, and hypergradient estimation.


\subsection{Synthetic Bilevel Problem}
Consider the bilevel optimization problem
\begin{align*}   
&\min\limits_{x\in \mathbb{R}^p} \Phi(x) := \frac{1}{|\mathcal{D}_{\text{val}}|}\sum\nolimits_{(u_i, v_i)\in \mathcal{D}_{\text{val}}}\frac{1}{2}\left((y^*(x))^\top u_i - v_i\right)^2 + (x^\top u_i - v_i)^2\\ \nonumber
&\text{s.t.} ~ y^*(x) :=  { \underset {y\in \mathbb{R}^p} { \operatorname {arg\,min} } \, g(x, y)} :=  \frac{1}{|\mathcal{D}_{\text{tr}}|}\sum\nolimits_{(u_i, v_i)\in \mathcal{D}_{\text{tr}}}\frac{1}{2}(y^\top u_i - v_i)^2 + \frac{r}{2}\| y - x \|^2
\end{align*}
where $r>0$, $\mathcal{D}_{\text{tr}}$ and $\mathcal{D}_{\text{val}}$ are datasets. For each given $x$, the LL function $g(x, \cdot)$ is strongly convex w.r.t. $y$. Furthermore, the explicit expressions of $y^*(x)$ can be obtained. This allow us to calculate the values of $\Phi(x)$ and $\|\nabla \Phi(x)\|$ during iterations, thus testing the numerical performance of algorithms.

In the experiments, we set $r$ to $0.5$, $p$ to $100$, and construct $\mathcal{D}_{\text{tr}}$ and $\mathcal{D}_{\text{val}}$, each with a sample size of $5000$, as follows: initially, we sample $10,000$ data points $e_i\in \mathbb{R}^{99}$, $i=1, \ldots, 10000$, from a normal distribution with a mean of $0$ and a variance of $0.01$. We then set $u_i = (e_i, 1)$, and $v_i = w_0^\top u_i + \sigma_i$ for each $i$, where $w_0 =(4, 6, 3, \ldots, 3)\in \mathbb{R}^{100}$, $\sigma_i \in \mathbb{R}$ represents the Gaussian noise with a mean of $0$ and a variance of $1$. The dataset $\{(u_i, v_i)\}_{i=1}^{10000}$ is then split equally into two parts to obtain $\mathcal{D}_{\text{tr}}$ and $\mathcal{D}_{\text{val}}$.

For VRBO, ALS-SPIDER, and ALS-STORM, we set the parameters as follows: $T=5$, $J=2$, $\alpha=0.01$, $\beta=0.1$, $\eta=0.01$, and $S_1=500$. To compare VRBO and ALS-SPIDER, we specify $q_1 = 10$ and $S_2 = 10$ for both algorithms. In addition, for the comparison between ALS-SPIDER and ALS-STORM, we let $q_1$ equal to $S_2$ for ALS-SPIDER, and set $\tau^x = 0.01$, $\tau^y = 0.0001$, and $\tau^v = 0.01$ for ALS-STORM.

\begin{figure}[htb]
  \centering
  \subfloat{\includegraphics[scale = 0.15]{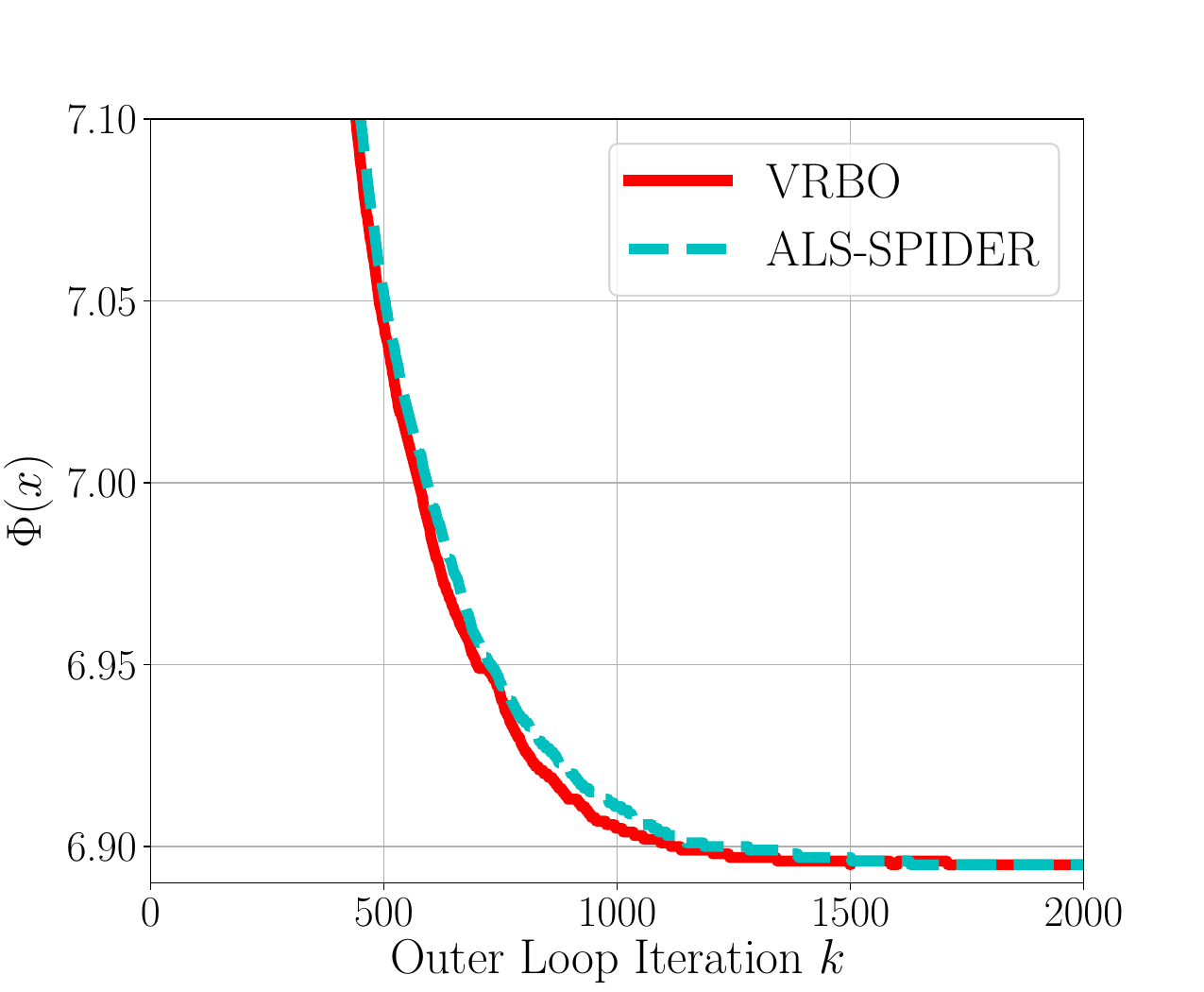}}%
  \hfil
   \subfloat{\includegraphics[scale = 0.15]{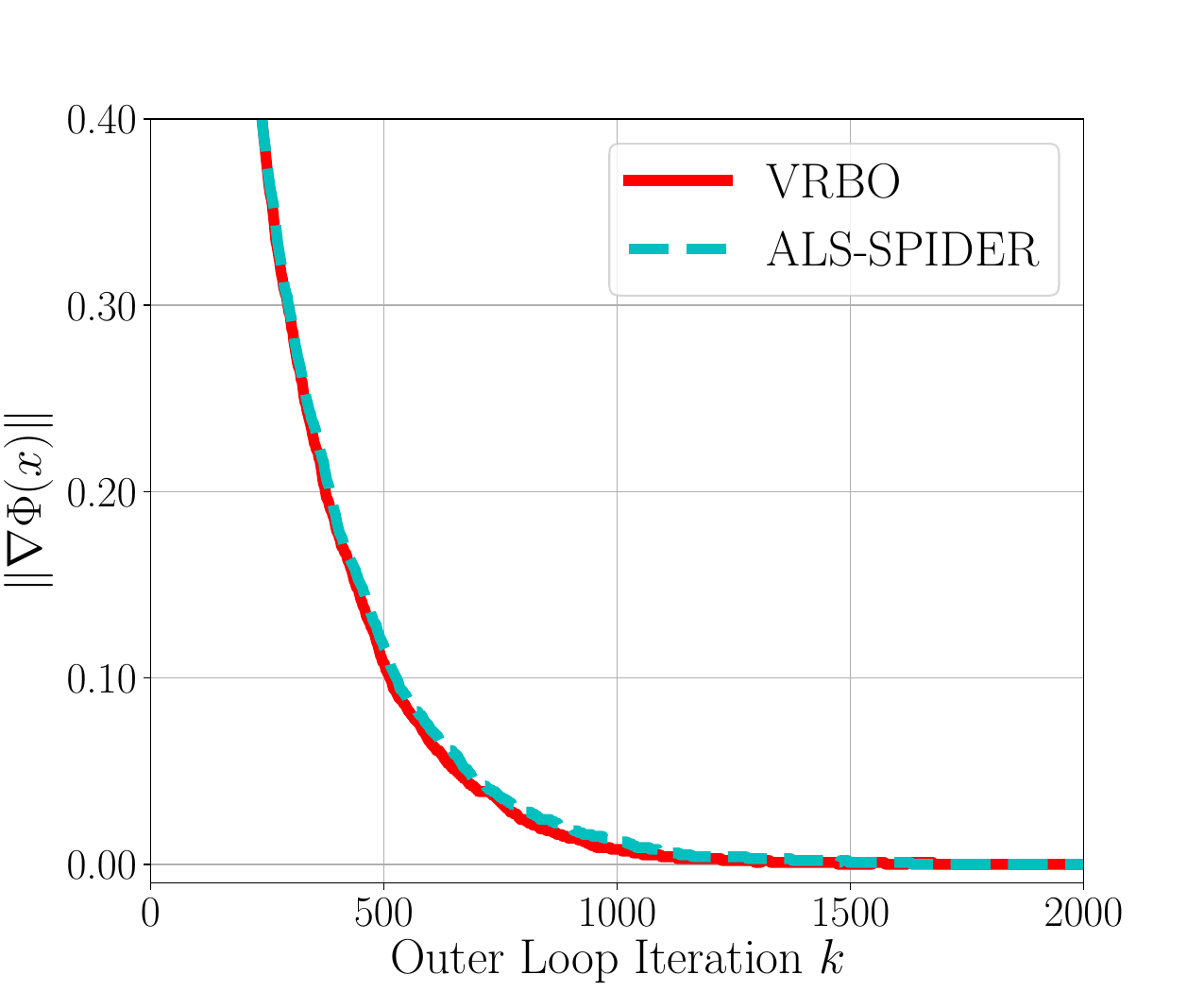}}%
  \hfil
   \subfloat{\includegraphics[scale = 0.15]{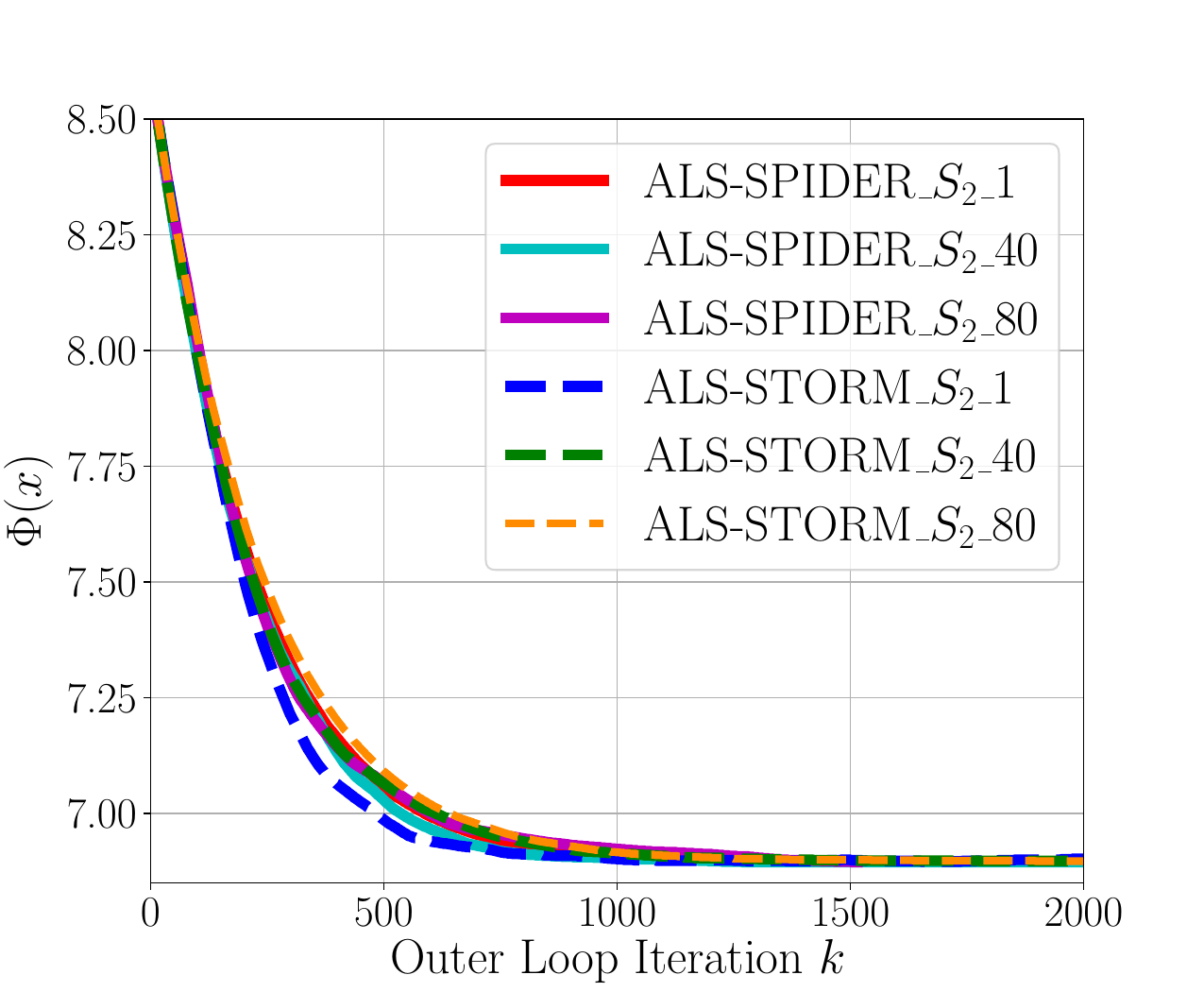}}%
    \hfil
   \subfloat{\includegraphics[scale = 0.15]{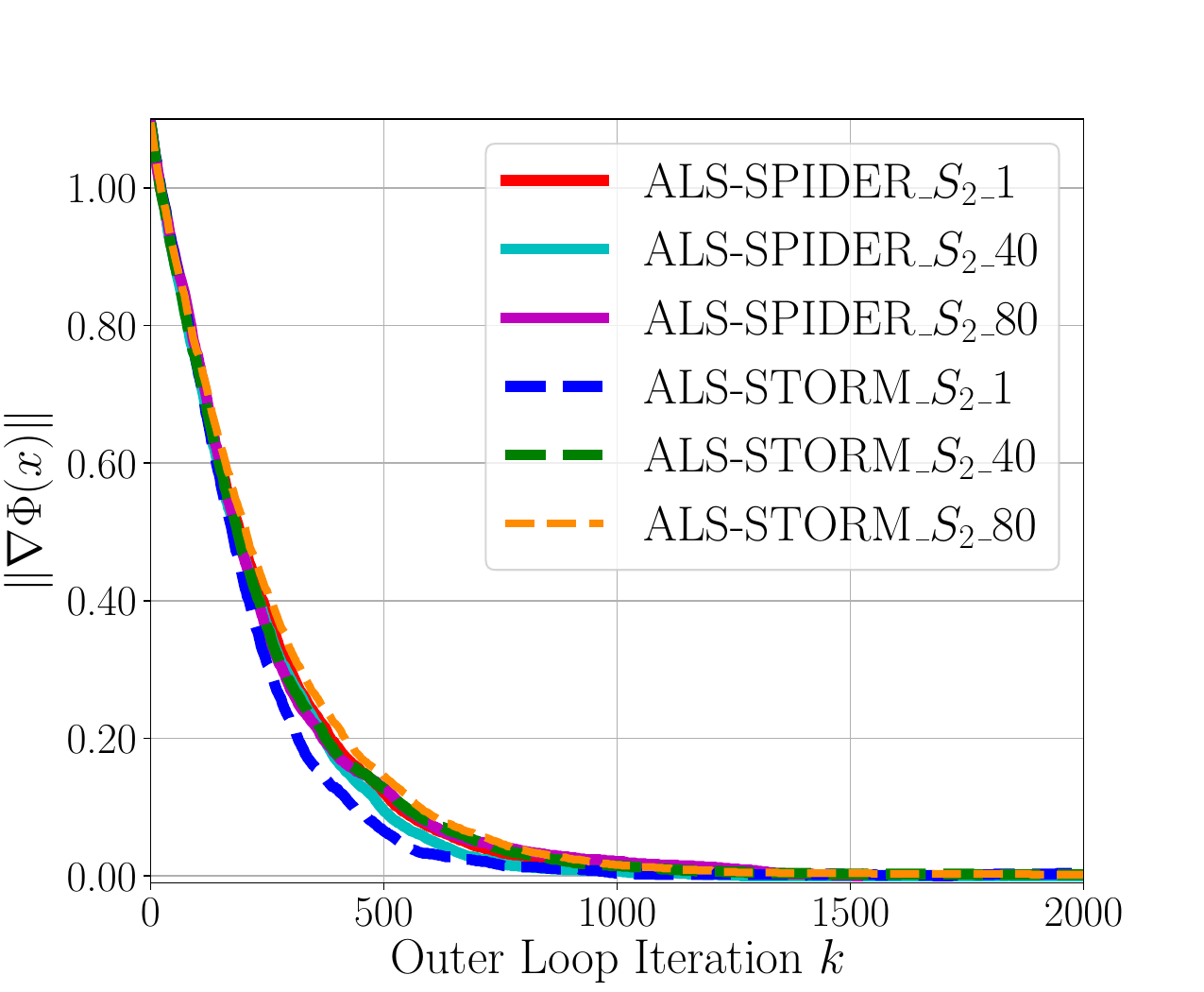}}%

   \subfloat{\includegraphics[scale = 0.15]{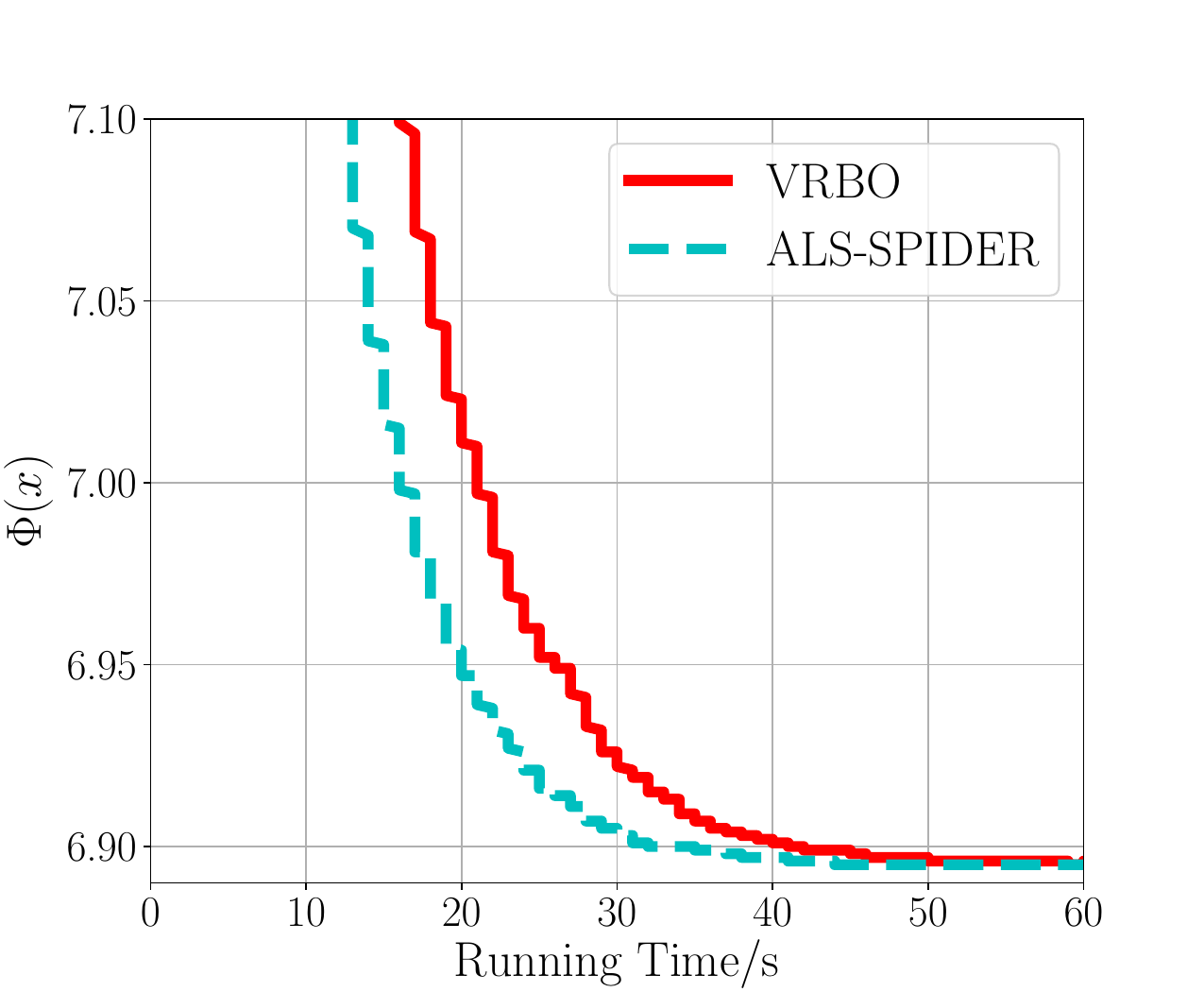}}%
  \hfil
   \subfloat{\includegraphics[scale = 0.15]{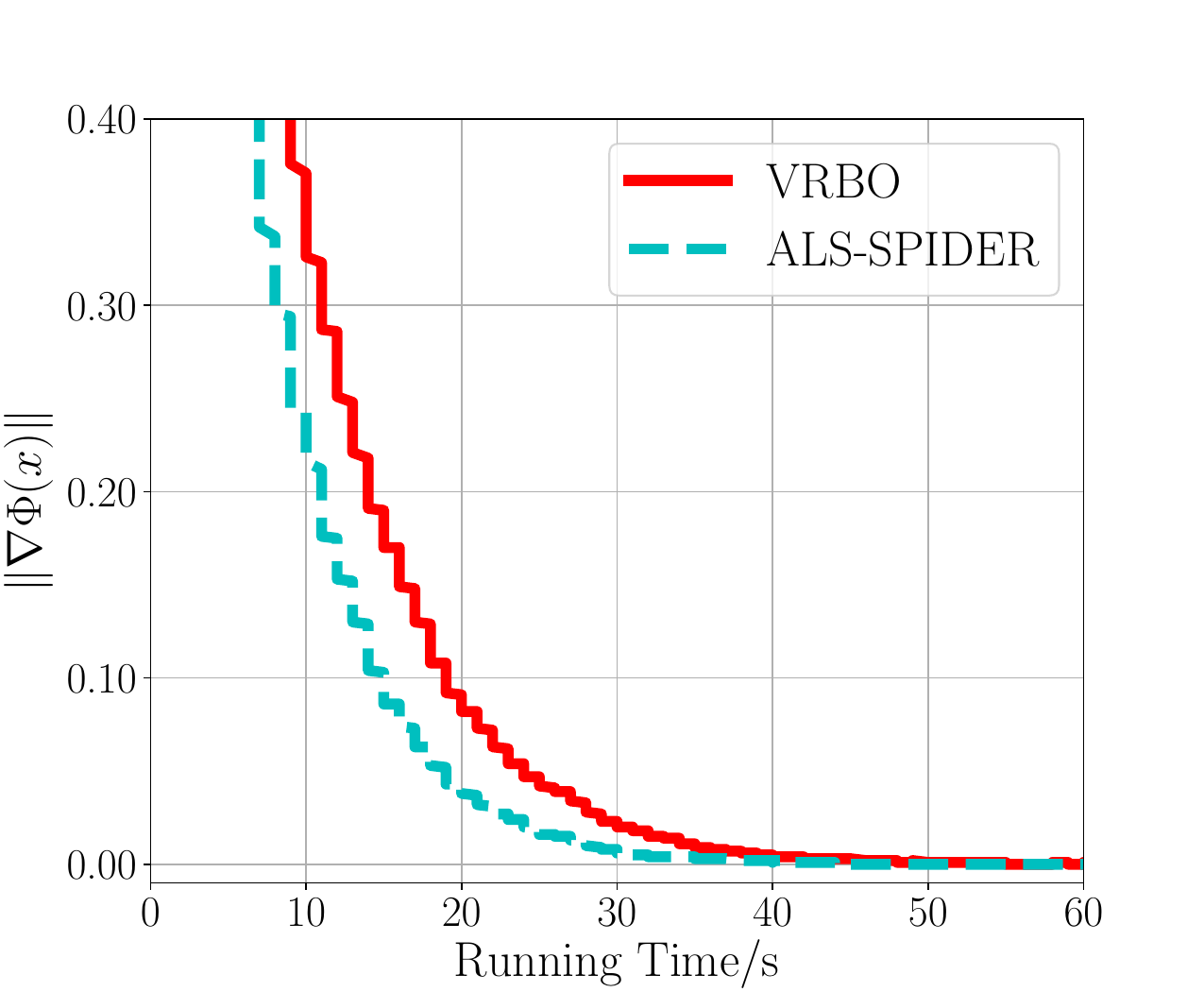}}%
  \hfil
   \subfloat{\includegraphics[scale = 0.15]{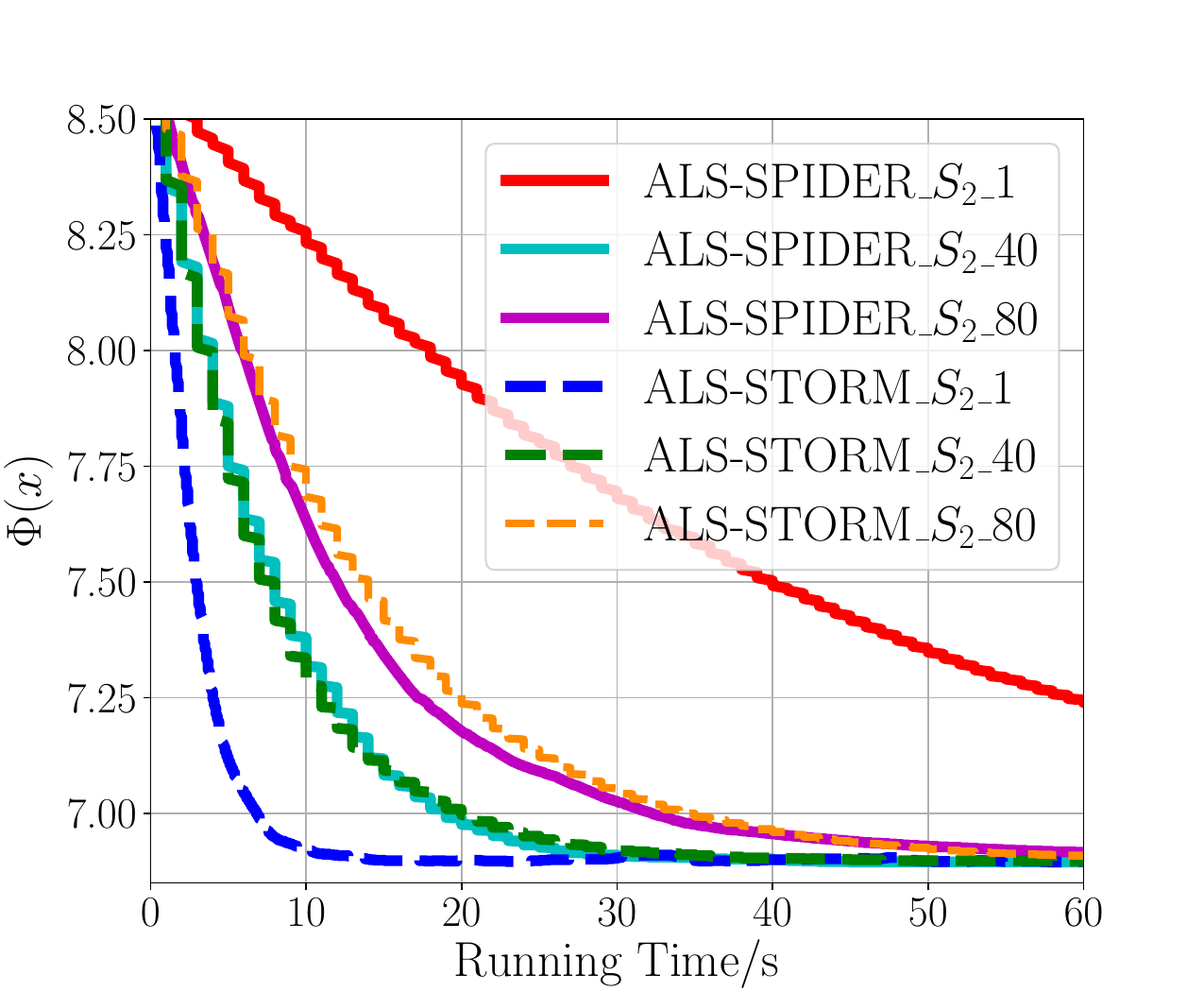}}%
    \hfil
   \subfloat{\includegraphics[scale = 0.15]{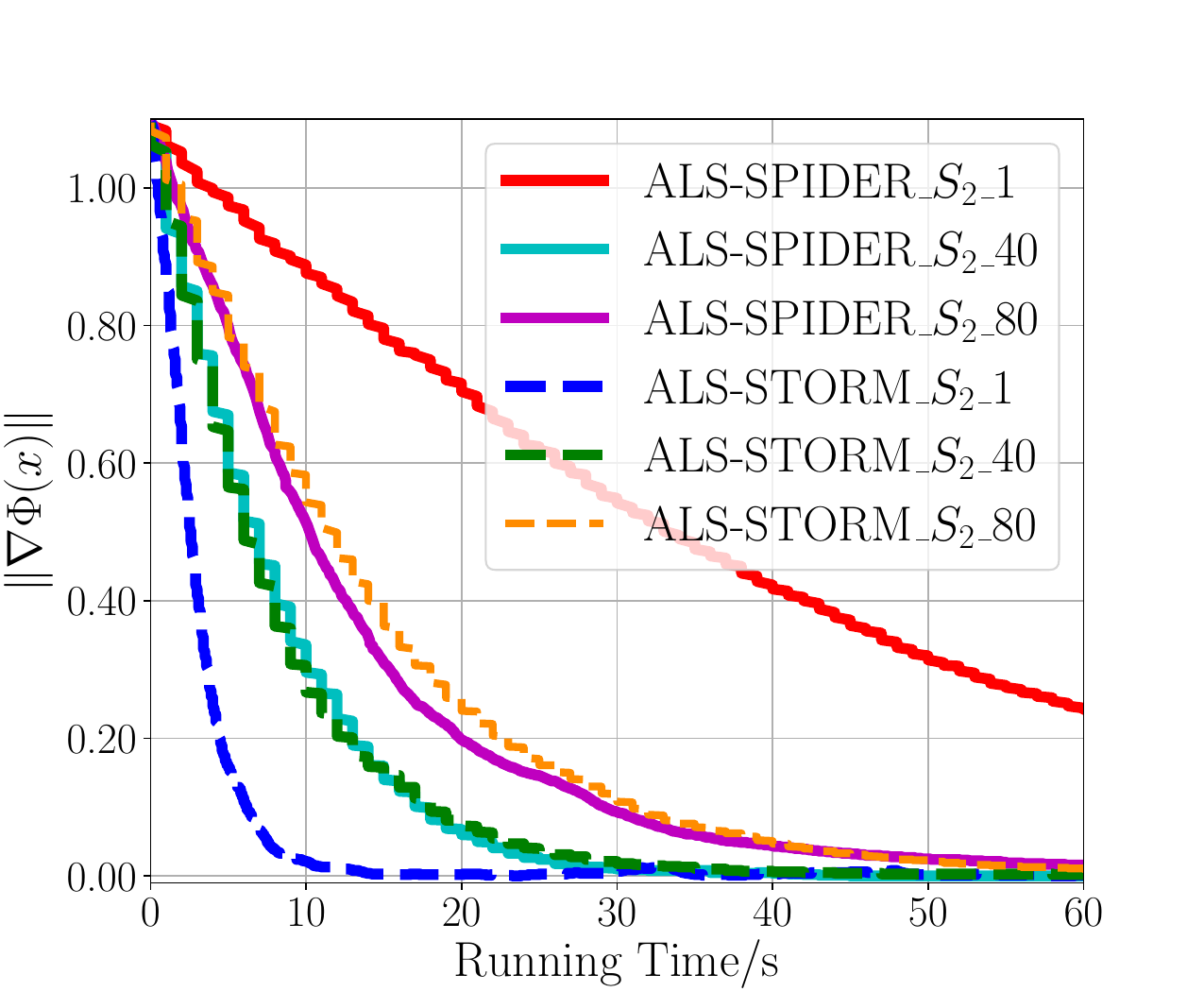}}%
  \caption{Comparison of VRBO and ALS-SPIDER (in the first two columns), and ALS-SPIDER and ALS-STORM (in the last two columns) for the synthetic bilevel problem. We show the convergence performance w.r.t. outer loop iteration $k$ (resp. running time) in the first (resp. second) row.}
  \label{fig:1}
\end{figure}

The comparisons between VRBO and ALS-SPIDER are presented in the first two columns of \cref{fig:1}. It can be observed that VRBO and ALS-SPIDER exhibit comparable performance in terms of the outer loop iteration $k$, w.r.t. both the function value $\Phi(x)$ and the gradient norm $\|\nabla \Phi(x)\|$. However, VRBO converges more slowly than ALS-SPIDER in terms of running time, primarily due to its complex three-level nested loop structure.

In the last two columns of \cref{fig:1}, we present the convergence results of ALS-SPIDER and ALS-STORM under three choices of $S_2 \in \{1, 40, 80\}$. Here, ALS-SPIDER\_$S_2$\_$i$ (resp. ALS-STORM\_$S_2$\_$i$) indicates that the batch size $S_2$ for ALS-SPIDER (resp. ALS-STORM) is set to $i$. It can be seen that even without using a periodic large batch size $S_1$, ALS-STORM performs comparably to ALS-SPIDER at each outer loop iteration $k$. Moreover, regarding running time, ALS-SPIDER achieves the fastest convergence rate with a smaller $S_2$ among $S_2 \in \{1, 40, 80\}$. Notably, with this smaller $S_2$, ALS-STORM achieves a faster convergence rate than ALS-SPIDER.

%

\subsection{Data Hyper-Cleaning}
Suppose there is a training set, where the label of each sample in the training set is corrupted with a probability $p$. In order to train a good classifier on this corrupted training set, data hyper-cleaning involves finding a suitable weight for each training sample to make the trained classifier minimize the loss on the validation set. When training a linear classifier, data hyper-cleaning can be expressed as a bilevel optimization problem in the following form\cite{13}
\begin{align} \label{eq75}
&\min\limits_{x} \Phi(x) := {1}/{|\mathcal{D}_{\text{val}}|}\sum\nolimits_{(u_i, v_i)\in \mathcal{D}_{\text{val}}}L\left(y^*(x); u_i, v_i\right) \\ \nonumber
&\text{s.t.} ~ y^*(x) :=  { \underset {y} { \operatorname {arg\,min} } \, g(x, y)} := {1}/{|\mathcal{D}_{\text{tr}}|}\sum\nolimits_{(u_i, v_i)\in \mathcal{D}_{\text{tr}}}\sigma(x_i)L(y; u_i, v_i) + c\|y\|^2
\end{align}
where $L$ denotes the cross-entropy loss, $\mathcal{D}_{\text{tr}}$ is the training set with a corruption probability $p$, $\mathcal{D}_{\text{val}}$ is the clean validation set, $\sigma(\cdot)$ is the sigmoid function, $y$ denotes the parameters of the classifier, $x :=\{x_i\}_{i=1}^{|\mathcal{D}_{\text{tr}}|}$, with $\sigma(x_i)$ being the weight for the training sample $(u_i, v_i)$, and $c$ is the coefficient of the regularization term.

We conduct experiments on the FashionMNIST\cite{27} dataset to solve problem (\ref{eq75}). In these experiments, we use $55,000$ images for training, $5,000$ for validation, $10,000$ for testing, with $c=0.01$, $p=0.3$. For all algorithms, we set $T = 5$, $J = 5$, $\alpha = 0.0001$, $\beta = 0.01$, $\eta = 0.001$, and $S_1 = 5000$. For the comparison between VRBO and ALS-SPIDER, we use $q_1 = S_2 = 5$ for both algorithms. In the comparison of ALS-SPIDER and ALS-STORM, $q_1$ is set equal to $S_2$ for ALS-SPIDER, while for ALS-STORM, we set $\tau^x = 0.01$, $\tau^y = 0.0001$, and $\tau^v = 0.01$.

\begin{figure}[htb]
  \centering
  \subfloat{\includegraphics[scale = 0.15]{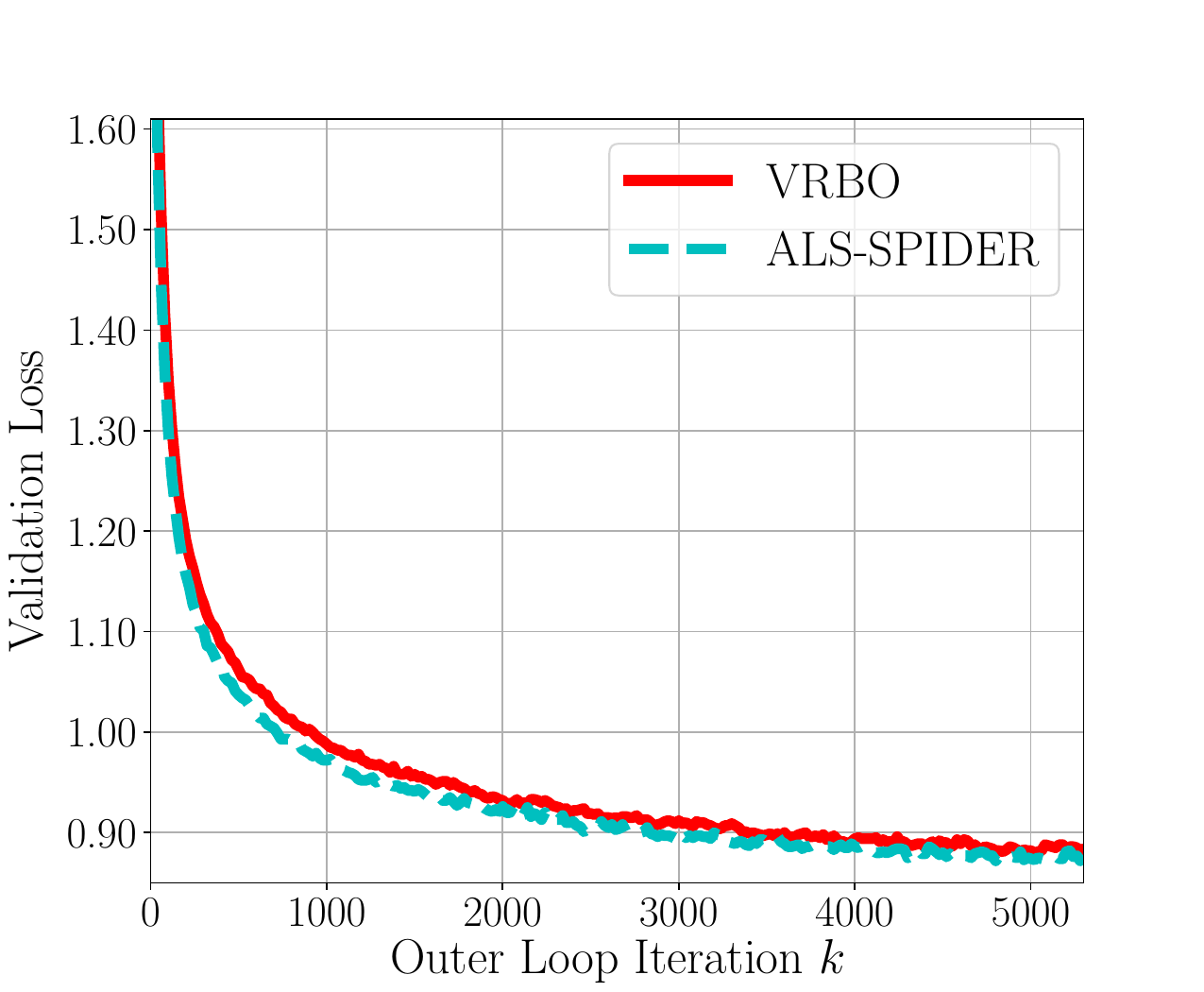}}%
  \hfil
   \subfloat{\includegraphics[scale = 0.15]{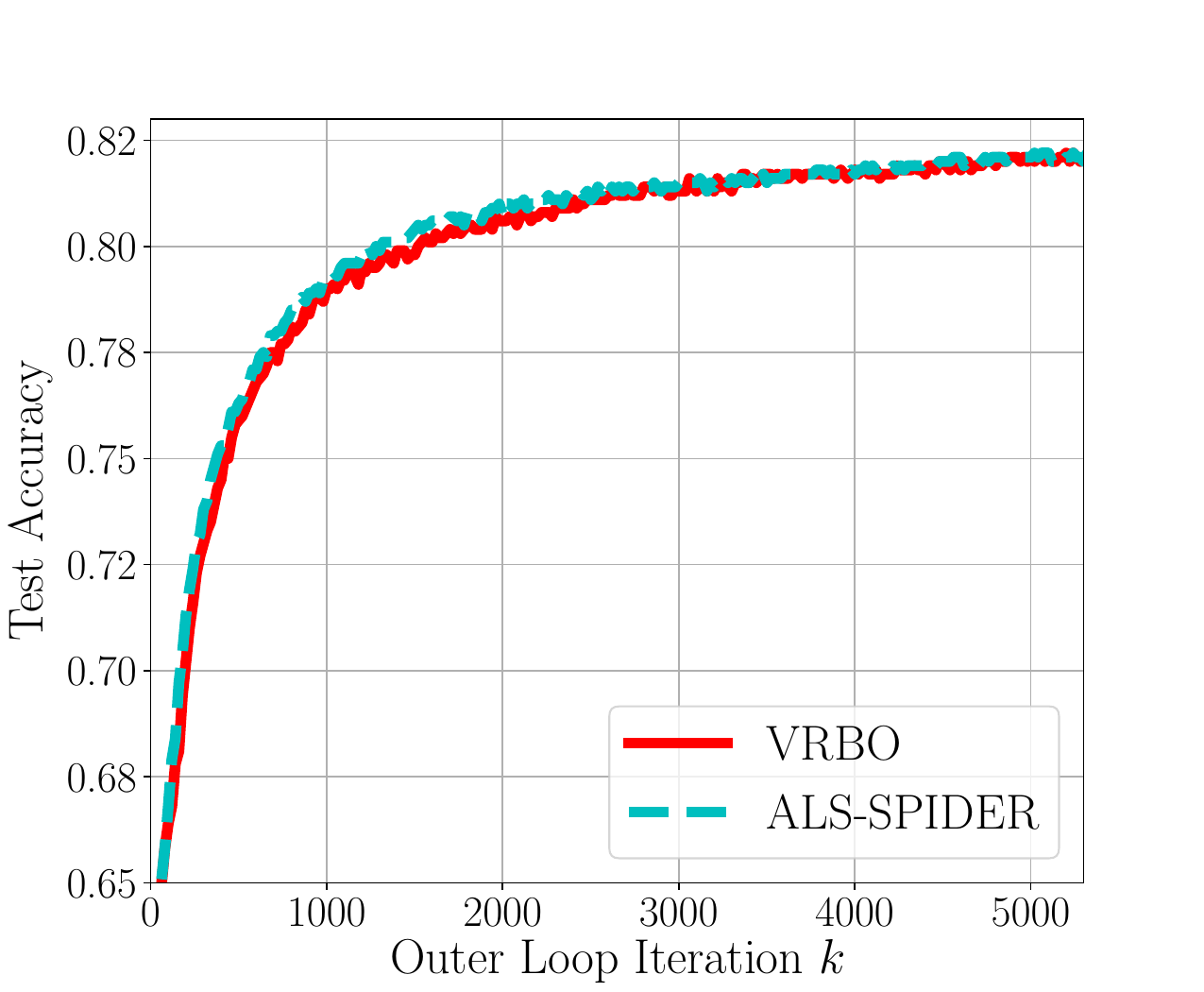}}%
  \hfil
   \subfloat{\includegraphics[scale = 0.15]{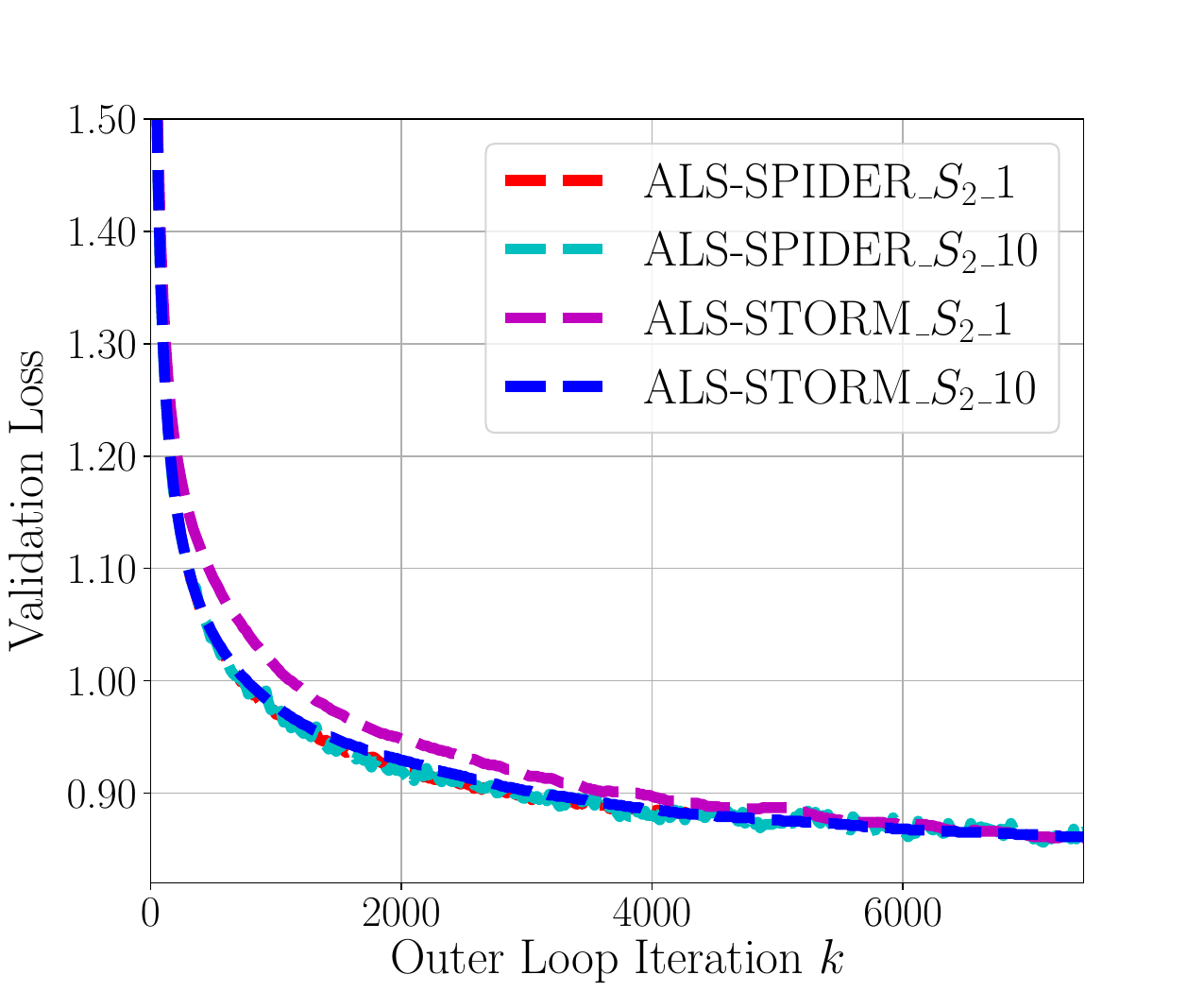}}%
    \hfil
   \subfloat{\includegraphics[scale = 0.15]{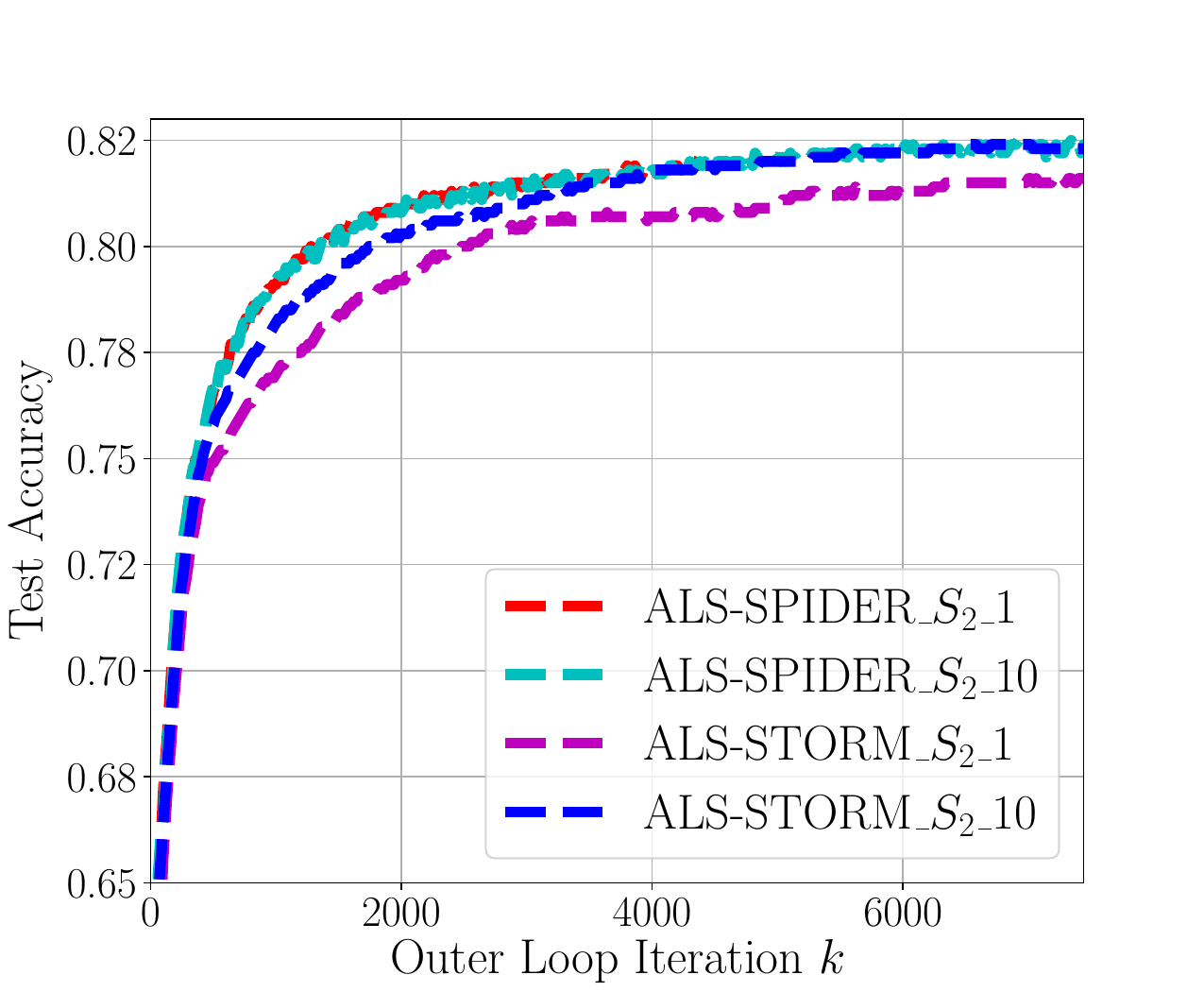}}%

   \subfloat{\includegraphics[scale = 0.15]{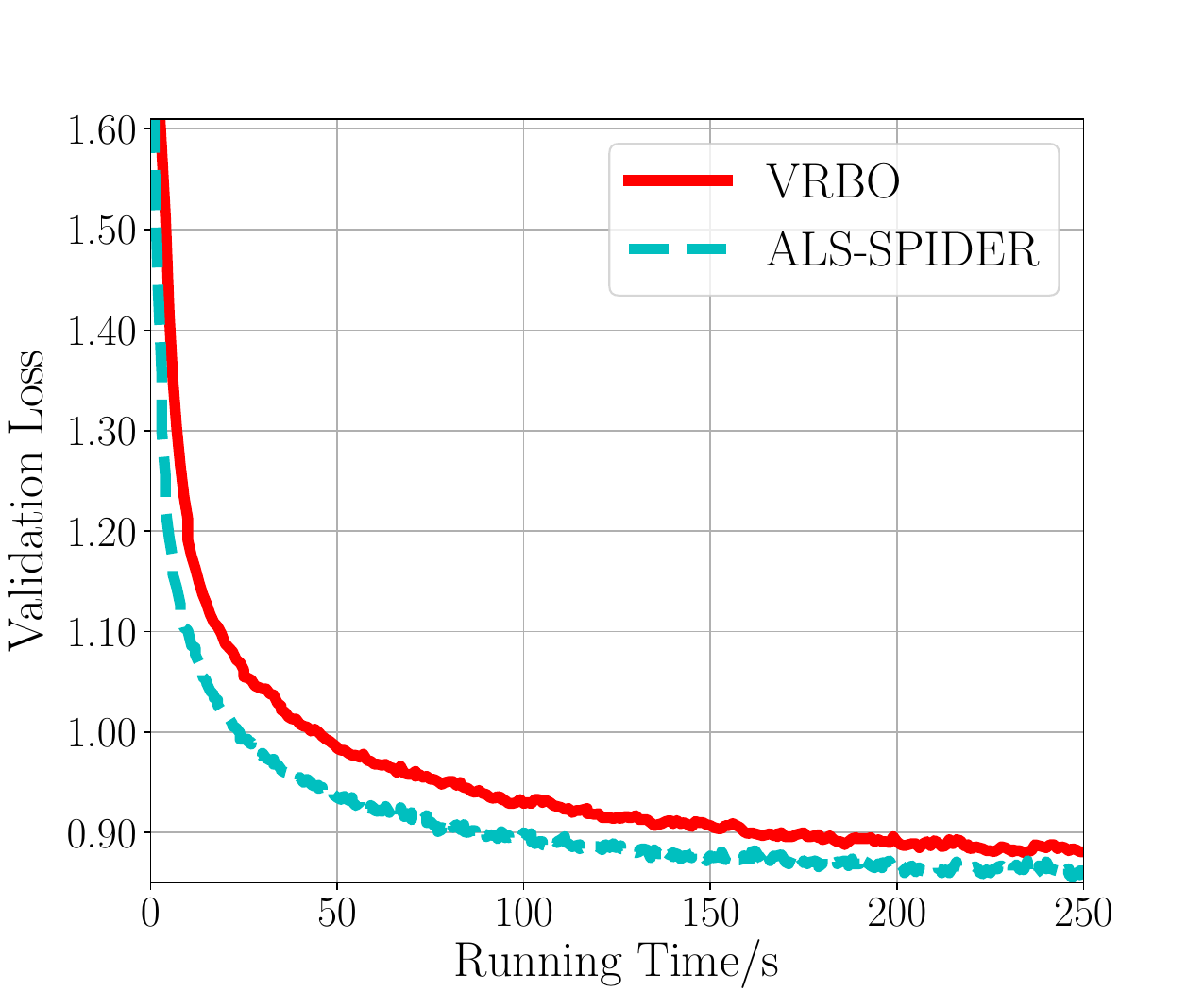}}%
  \hfil
   \subfloat{\includegraphics[scale = 0.15]{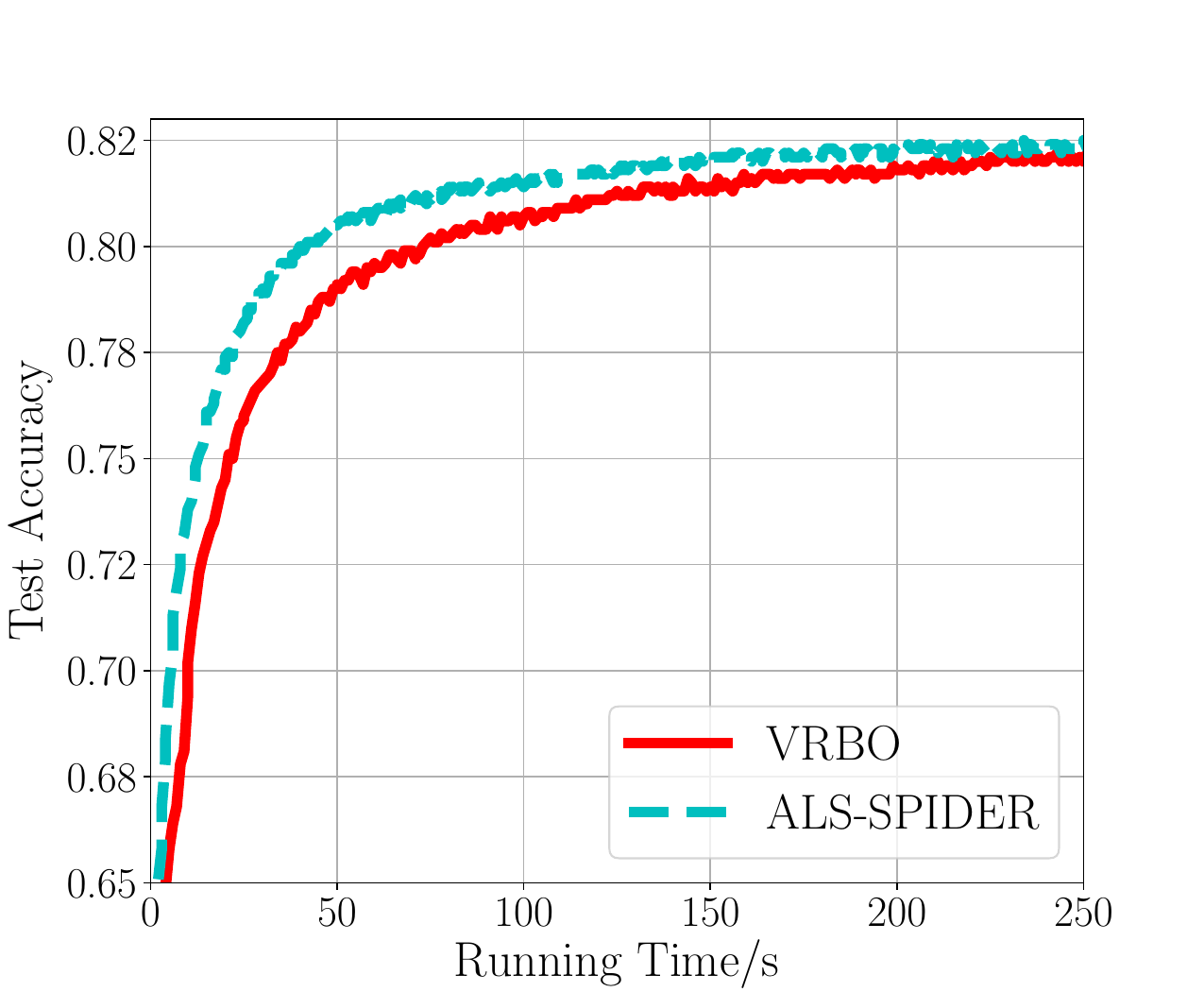}}%
  \hfil
   \subfloat{\includegraphics[scale = 0.15]{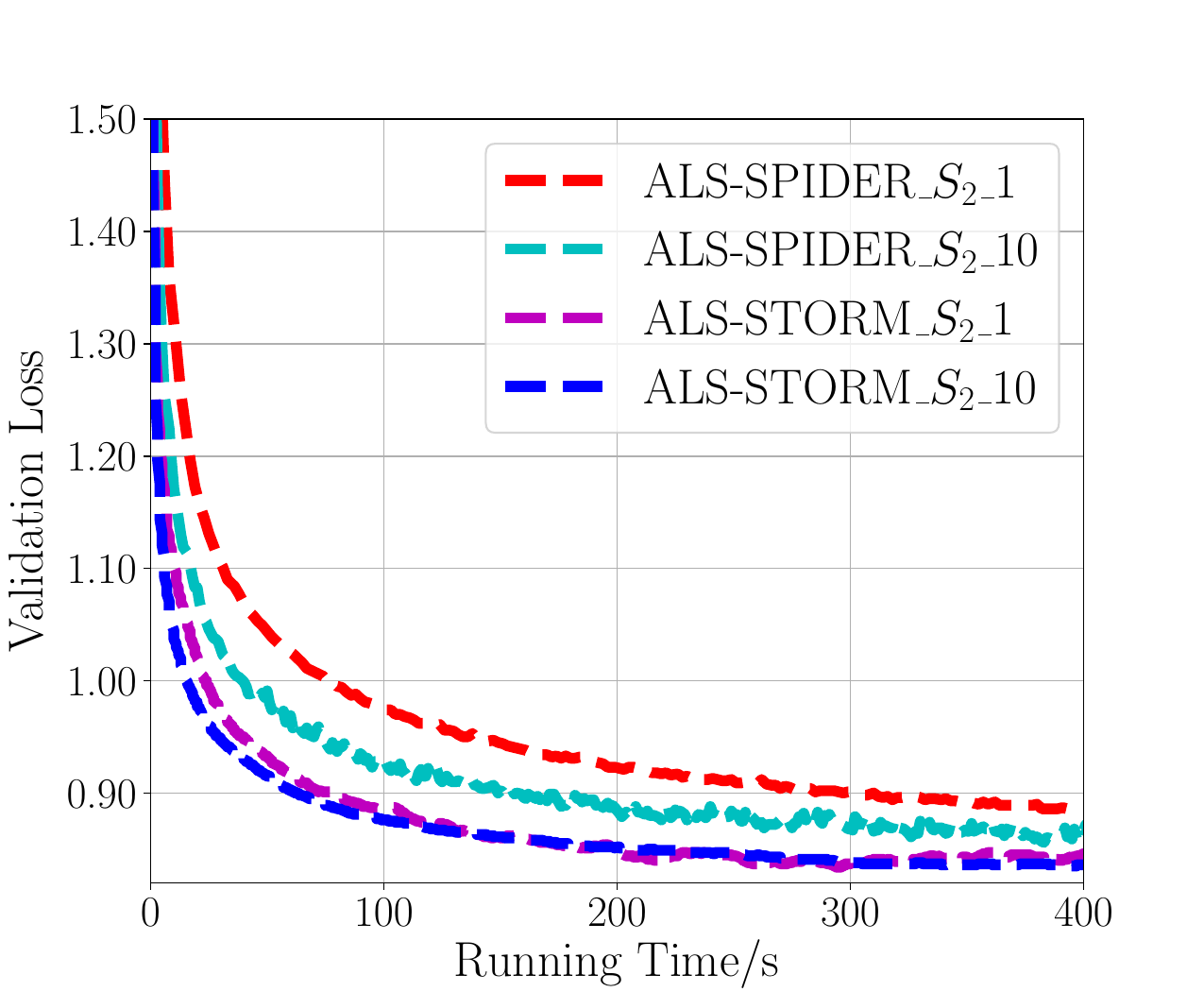}}%
    \hfil
   \subfloat{\includegraphics[scale = 0.15]{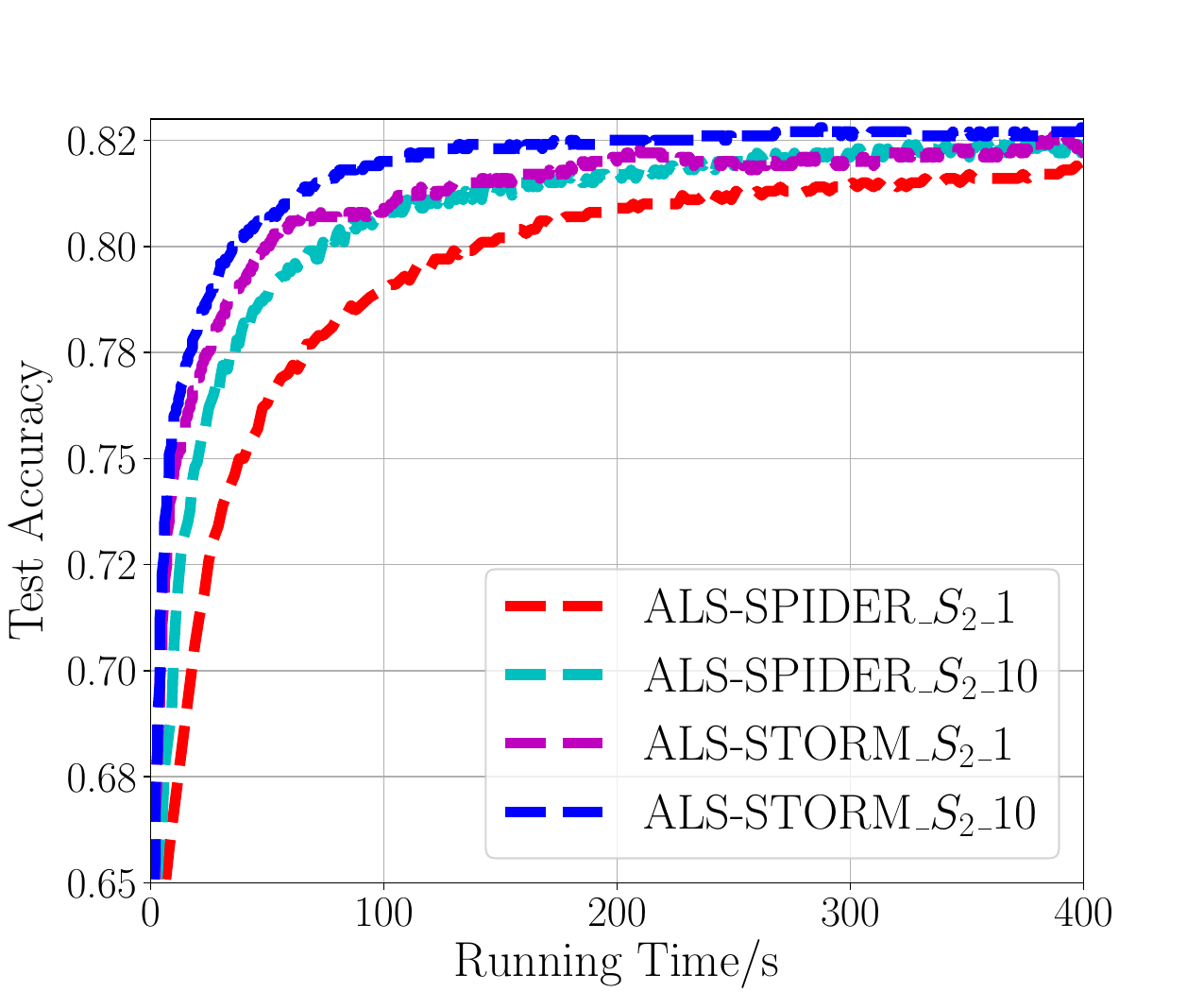}}%
  \caption{Comparison of VRBO and ALS-SPIDER (in the first two columns), and ALS-SPIDER and ALS-STORM (in the last two columns) for the data hyper-cleaning task on the FashionMNIST dataset. Results w.r.t. outer iteration $k$ (resp. running time) are in the first (resp. second) row. }
  \label{fig:2}
\end{figure}

The experimental results are shown in \cref{fig:2}, where the validation loss on $\mathcal{D}_{\text{val}}$, and the test accuracy on the test set are calculated based on the estimated $y$ from the training set $\mathcal{D}_{\text{tr}}$, and ALS-SPIDER\_$S_2$\_$i$ and ALS-STORM\_$S_2$\_$i$ are definined as before. Consistent with the results for the synthetic bilevel problem, ALS-SPIDER outperforms VRBO in terms of running time. In addition, it is observed that a larger $S_2$ accelerates the convergence of ALS-SPIDER w.r.t. running time, due to the lower sample complexity. In constrast, the convergence rate of ALS-STORM w.r.t. running time is less sensitive to the value of $S_2$. Notably, when $S_2 = 10$, ALS-STROM achieves the fastest convergence rate, the lowest validation loss, and the highest test accuracy among the results in the last two columns of \cref{fig:2}, in terms of running time. This highlights the practical advantange of ALS-STORM over ALS-SPIDER.

\section{Conclusion}\label{sec:6}
In this paper, we propose two alternating stochastic variance-reduced algorithms, ALS-SPIDER and ALS-STORM, for unconstrained nonconvex-strongly-convex bilevel optimization, by introducing an auxiliary variable for hypergradient estimation. ALS-SPIDER employs the SPIDER variance reduction technique, whereas ALS-STORM is a variant of ALS-SPIDER to avoid using  large batches in every iteration. Both algorithms achieve the optimal complexity of $O(\epsilon^{-1.5})$ to reach an $\epsilon$-stationary point when updating the lower-level variable multiple times, outperforming existing alternating stochastic algorithms with multi-step lower-level variable updates. Experimental results showcase their superior performance.

\bibliographystyle{unsrt}
\bibliography{mybible}

\end{document}